\documentclass[a4paper,11pt]{article}
\usepackage[latin1]{inputenc}
\usepackage{amssymb}
\usepackage{amsmath}
\usepackage[english]{babel}
\usepackage[colorlinks, linkcolor=red, urlcolor=blue, citecolor=blue]{hyperref}
\newtheorem{thm}{Theorem}[section]

\newtheorem{coro}{Corollary}[section]
\newtheorem{prop}{Proposition}[section]
\newtheorem{defn}{Definition}[section]
\newtheorem{exple}{Example}[section]
\newtheorem{rem}{Remark}[section]
\newenvironment{proof}[1][Proof]{\textbf{#1.}}{\hfill $\Box$}
\begin{document}

\date{}
\title{\LARGE{IDT processes and associated Lévy processes with explicit constructions}\footnote{This work is partially supported by Hassan II Academy of Sciences and Technology.}}
\author{Antoine HAKASSOU\footnote{\textbf{a.hakassou@edu.uca.ma}, Cadi Ayyad University, LIBMA Laboratory, Department of Mathematics,
Faculty of Sciences Semlalia P.B.O. 2390 Marrakesh, Morocco.} \, and \, Youssef OUKNINE\footnote{\textbf{ouknine@uca.ma}, Cadi Ayyad University of
Marrakesh and Hassan II Academy of Sciences and Technology Rabat.}
}

\maketitle

\begin{center}
 \textbf{Abstract}
\end{center}

This article deals with IDT processes, i.e. processes which are infinitely divisible with respect to time. 
Given an IDT process $(X_{t},\,t\geq0)$, there exists a unique (in law) Lévy process $(L_{t}; t\geq0)$ which has the same one-dimensional marginals 
distributions, i.e. for any $t\geq0$ fixed, we have $$X_{t}\stackrel{(law)}{=}L_{t}.$$
Such processes are said to be associated. The main objective of this work is to exhibit numerous examples of IDT processes and their associated Lévy 
processes. To this end, we take up ideas of the monograph \textit{Peacocks and associated martingales} from F. Hirsch, C. Profeta, B. Roynette and M. 
Yor (Lévy, Sato and Gaussian sheet methods) and apply them in the framework of IDT processes. This gives a new interesting outlook to the study of 
processes whose only one-dimensional marginals are known. Also, we give an integrated weak Itô type formula for IDT processes 
(in the same spirit as the one for Gaussian processes) and some links between IDT processes and selfdecomposability. 
The last sections are devoted to the study of some extensions of the notion of IDT processes in the weak sense as well as in the multiparameter sense.
In particular, a new approach for multiparameter IDT processes is introduced and studied. 
Main examples of this kind of processes are the $\mathbb{R}_{+}^{N}-$parameter L\'{e}vy process and 
the L\'{e}vy's $\mathbb{R}^{M}$-parameter Brownian motion. These results give some 
better understanding of IDT processes, and may be seen as some continuation of the works of K. Es-Sebaiy and Y. Ouknine 
[\textit{How rich is the class of processes which are infinitely divisible with respect to time ?}] and R. Mansuy 
[\textit{On processes which are infinitely divisible with respect to time}]. \\



\section{Introduction}

In this introductory section, we are going to give several examples of IDT processes (essentially constructed from Gaussian processes), as well as 
associated Lévy processes. First, let us recall some well known facts on the class of IDT processes.

\begin{defn}

A stochastic process $X=(X_{t}; t\geq 0)$ is an IDT process if it satisfies the following condition:
\begin{equation}\label{1}
\forall n \in\mathbb{N}^{*}, (X_{nt}; t\geq 0) {\stackrel{(law)}{=}}(X_{t}^{(1)}+\cdots+X_{t}^{(n)}; t\geq 0)
\end{equation}
where $(X^{(i)})_{1\leq i\leq n}$ are independent copies of X.

\end{defn}

\begin{exple}

Let $0<\alpha<2$ and consider an $(\alpha/2)$-stable positive random variable $\xi$. Now, we consider a centered Gaussian process $(G_{t};t\geq0)$, 
such that its covariance function $R(s,t):=\mathbb{E}[G_{s}G_{t}]$ verifies for all $\lambda>0$, $R(\lambda s, \lambda t)=\lambda^{2/\alpha}R(s,t)$. 
Assume that $(G_{t};t\geq0)$ is stochastically continuous and independent of $\xi$. \\ Then, the stochastic process $(X_{t};t\geq0)$ defined 
by $$X_{t}:=\xi^{1/2}G_{t},$$ is an IDT process, and following P. Embrechts and M. Maejima (\cite{embrechts2002selfsimilar} Example 3.6.4), it is also a 
stable sub-Gaussian process.\par
In fact, since $\xi$ is an $(\alpha/2)$-stable variable, we have $\mathbb{E}[exp\{-z\xi\}]=exp\{-|z|^{\alpha/2}\}$ for any $z\in\mathbb{R}$. Then, for all $n,m\in\mathbb{N}^{\ast}$, we get:
$$\mathbb{E}[exp\{i \theta\sum_{k=1}^{m} a_{k} X_{n t_{k}}\}]=\mathbb{E}[exp\{i \theta\sum_{k=1}^{m} a_{k} \xi^{1/2} G_{n t_{k}}\}].$$
This implies that:
$$\mathbb{E}[exp\{i \theta \sum_{k=1}^{m} a_{k} X_{n t_{k}}\}]=\mathbb{E}_{\xi} \mathbb{E}_{G} [exp\{i \theta \xi^{1/2} \sum_{k=1}^{m} a_{k} G_{n t_{k}}\}].$$
Then, we get:
$$\mathbb{E}[exp\{i \theta \sum_{k=1}^{m} a_{k} X_{n t_{k}}\}]=\mathbb{E}_{\xi} [exp \{-\frac{1}{2} {|\theta|}^{2} \xi  \sum_{k,j=1}^{m} a_{k} a_{j} R(n t_{k}, n t_{j})\}].$$
And then, 
$$\mathbb{E}[exp\{i \theta \sum_{k=1}^{m} a_{k} X_{n t_{k}}\}]=exp \{-{|\theta|}^{\alpha} [\frac{1}{2}\sum_{k,j=1}^{m} a_{k} a_{j} R(n t_{k}, n t_{j})]^{\alpha/2}\}.$$
According to the properties of the covariance function R, we get:
$$\mathbb{E}[exp\{i \theta \sum_{k=1}^{m} a_{k} X_{n t_{k}}\}]=exp \{-n{|\theta|}^{\alpha} [\frac{1}{2}\sum_{k,j=1}^{m} a_{k} a_{j} R(t_{k},t_{j})]^{\alpha/2}\}.$$
Hence,
$$\mathbb{E}[exp\{i \theta \sum_{k=1}^{m} a_{k} X_{n t_{k}}\}]=(\mathbb{E}[exp\{i \theta \sum_{k=1}^{m} a_{k} X_{t_{k}}\}])^{n}.$$

\end{exple}

Following R. Mansuy \cite{Mansuy}, one may wonder which among centered Gaussian processes $(G_{t}; t\geq0)$ (which, for simplicity, are assumed to be centered) are IDT.
In order to give a characterization of IDT Gaussian processes, we recall the following proposition due to R. Mansuy (\cite{Mansuy} Proposition 3.2).

\begin{prop}

Let $(G_{t}; t\geq0)$ be a centered Gaussian process, which is assumed to be continuous in probability. Then the following properties are equivalent:
\begin{enumerate}
\item  $(G_{t}; t\geq 0)$ is an IDT process.
\item  The covariance function $c(s,t):=\mathbb{E}[G_{s}G_{t}]$, $0 \leq s \leq t,$ satisfies
                   $$\forall \alpha > 0,  c(\alpha s, \alpha t)=\alpha c(s,t), \mbox{ for all   }  0 \leq s \leq t.$$
\item  The process $(G_{t}; t\geq 0)$ satisfies the "Brownian scaling property", namely
          $$\forall \alpha > 0, \mbox{  } (G_{\alpha t}; t \geq 0)\stackrel{(law)}{=}(\sqrt{\alpha} G_{t}; t\geq 0)$$
\item  The process $(\tilde{G}_{y}:=e^{-y/2} G_{e^{y}}; y\in\mathbb{R})$ is stationary.
\item  The covariance function $\tilde{c}(y,z):=\mathbb{E}[\tilde{G}_{y} \tilde{G}_{z}],$  $y,z\in\mathbb{R},$ is of
the form $$\tilde{c}(y,z)=\int \mu (du) e^{iu \mid y-z \mid},  \mbox{ y,z}\in\mathbb{R}$$
where $\mu$ is a positive, finite, symmetric measure on $\mathbb{R}$.
\end{enumerate}
Then, under these equivalent conditions, the covariance function c of $(G_{t}; t\geq 0)$ is given by $$c(s,t)=\sqrt{st} \int \mu (da) e^{ia \mid ln(\frac{s}{t}) \mid}.$$

\end{prop}

\begin{exple}

Let $(G_{t};t\geq 0)$ be a centered Gaussian process such that its covariance function $c(s,t):=\mathbb{E}[G_{s}G_{t}]$, $0 \leq s \leq t,$
is continuous and satisfies for $0\leq\alpha\leq1$,
$$\forall \lambda > 0, c(\lambda t, \lambda s)= {\lambda}^{\alpha}c(s,t) \mbox{ for all } 0\leq s \leq t.$$
Then, the stochastic process $(\tilde{G}_{t}; t\geq 0)$ defined by
$$\tilde{G}_{t}:=t^{\frac{1-\alpha}{2}}G_{t}  \mbox{     for all  } t\geq 0$$
is an IDT Gaussian process. \\

\end{exple}

Now, in the following we are going to recall some relationships that bind IDT processes to Lévy processes.

\begin{thm}

Any Lévy process is an IDT process and conversely, any stochastically continuous IDT process with independent increments, is a Lévy process.

\end{thm}

\begin{proof}

The first implication is easy. For the second, it is enough to prove the stationary increments property. More details could be
found in K. Es-Sebaiy and Y. Ouknine (\cite{Ouknine} Theorem 3.1) and R. Mansuy (\cite{Mansuy} Proposition 1.1).

\end{proof}

\begin{defn}

A stochastic process $(X_{t}; t\geq 0)$ is associated to a stochastic process $(Y_{t}; t\geq 0)$ if they have the same one-dimensional marginals
distributions, i.e. $\mbox{ for any fixed } t\geq0, \,\,\, X_{t}\stackrel{(law)}{=}Y_{t}.$

\end{defn}

\begin{prop}

If X is a stochastically continuous IDT process, then there exists a unique (in law) associated Lévy process L, i.e. for any fixed $t\geq0,$ $X_{t}\stackrel{(law)}{=}L_{t}$.

\end{prop}

\begin{proof}

Given an IDT process $X=(X_{t}; t\geq 0)$, notice that $X_{1}$ is an infinitely divisible random variable.
Then there exist a unique (in law) Lévy process $L=(L_{t}; t\geq 0)$ such that $X_{1}\stackrel{(law)}{=}L_{1}.$
Hence, according to the stochastically continuity of X and L,  $X_{t}\stackrel{(law)}{=}L_{t}$ for any fixed $t\geq 0$.

\end{proof}

In what follows, we illustrate the above proposition and give several examples of IDT processes as well as associated Lévy processes respectively.

\begin{exple}

Let $\{X_{t};t\geq0\}$ be a non-trivial strictly $\alpha$-stable Lévy process on $\mathbb{R}$ with $0<\alpha<2$. 
Define $Y_{t}=t^{2/\alpha}X_{1/t}$ for $t>0$ and $Y_{0}=0$. Then, we show that $\{Y_{t};t\geq0\}$ is an IDT process which is associated to the Lévy process X.

\end{exple}

\begin{exple}

For $\alpha \leq \frac{1}{2}$, we consider the centered Gaussian process G defined by
$$(G_{t}; t\geq 0):=(t^{\alpha} B_{t^{1-2\alpha}}; t\geq 0)$$
where $(B_{t};t\geq0)$ is the standard Brownian motion.
Then, one can easily point that $$\mathbb{E}[G_{\lambda t}G_{\lambda s}]=\lambda \mathbb{E}[G_{t}G_{s}] \mbox{ for all } \lambda >0.$$
Hence $(G_{t};t\geq0)$ is an IDT Gaussian process and its associated Lévy process is the standard Brownian motion $(B_{t};t\geq0)$
since $\forall t\geq0,$ $Var(G_{t})=Var(B_{t})=t$.

\end{exple}

\begin{exple}

Let $(B^{H}_{t};t\geq0)$ be a standard fractional Brownian motion with Hurst index $H\in (0,1)$.
One can easily show that the stochastic process $\bar{G}$ defined by $$(\bar{G}_{t}; t\geq 0):=(t^{\frac{1}{2}-H}B^{H}_{t}; t\geq 0),$$ is an IDT Gaussian process
and its associated Lévy process is the standard Brownian motion $(B_{t};t\geq0)$.

\end{exple}

\begin{exple}

We consider the continuous Gaussian semimartingale $(X_{t}; t\geq 0)$ defined by $$X_{t}=B_{t}-\int_{0}^{t} \int_{0}^{u} l(u,v)d B_{v} du$$
where $(B_{t}; t\geq 0)$ is the standard Brownian motion, l a continuous Volterra kernel of the form $l(u,v)=\frac{1}{u} \varphi (\frac{v}{u})$ and $\varphi$
a function which satisfies
\begin{equation}\label{2}
\int_{0}^{1} \varphi (x)dx= \int_{0}^{1} \int_{0}^{1} \varphi (zx) \varphi (z) dz dx
\end{equation}
without satisfying
\begin{equation}\label{3}
\varphi (x) = \int_{0}^{1} \varphi (zx) \varphi (z) dz.
\end{equation}
An example of such functions $\varphi$ could be found in H. Föllmer, C.T. Wu and M. Yor (\cite{Follmer} Section 6.2 Theorem 6.3). To be quite explicit, they consider
$\varphi (x)=ce^{-ax}$ and then, for a given real number $a$, \eqref{2} is satisfied if and only if $$c=\frac{(1-e^{-a})}{\int_{0}^{a} e^{-u}(1-e^{-u})\frac{du}{u}},$$
whereas \eqref{3} is never satisfied for any $c\neq0$. \\
We further assume that
$$\int_{0}^{t}(\int_{0}^{u} l^{2}(u,v)dv)^{1/2} du=\int_{0}^{t}(\int_{0}^{u} \frac{1}{u^{2}} {\varphi}^{2}(\frac{v}{u})dv)^{1/2} du <+\infty.$$
Clearly X has quadratic variation $<X>_{t}=t$, and according to H. Föllmer, C.T. Wu and M. Yor (\cite{Follmer} Section 6.2.), $(X_{t};t\geq 0)$ is a weak
Brownian motion of order 1 which is not a Brownian motion. Now, one can easily show that $(X_{t};t\geq 0)$ is an IDT Gaussian process and its associated
Lévy process is the standard Brownian motion $(B_{t};t\geq 0)$.

\end{exple}

\begin{rem}

Consider a stochastically continuous IDT process X and L its associated Lévy process. Denote by $(b,\sigma,\nu)$ the characteristic triplet of L.
Then we can write
\begin{equation*}
L_{t}=bt+\sigma W_{t}+\int_{0}^{t} \int_{\mid x \mid \geq 1} x {\mu}^{L} (ds,dx)+
\int_{0}^{t} \int_{\mid x \mid < 1} x ({\mu}^{L} (ds,dx)-ds{\nu}(dx))
\end{equation*}
where ${\mu}^{L}(ds,dx)$ denote the random measure counting the jumps of L.\\
If $\mathbb{E}(X_{1})=\mathbb{E}(L_{1})=b+\int_{\mid x \mid \geq 1} x {\nu} (dx)=0$ and $\mathbb{E}|L_{1}| < +\infty$, then L is a martingale and X is a 1-martingale
and so X is a Peacock in the sense of F. Hirsch \textit{et al} \cite{Hirschbook}.

\end{rem}

In the following, we give a result which allows us to derive a Fokker-Planck PDE for IDT processes.

\begin{thm}

Consider X a stochastically continuous IDT process and denote by $p_{t}^{X}$ the law of $X_{t}$. Let L be the Lévy process (in law) associated to X and denote
by $\mathcal{L}^{*}$ the adjoint of its infinitesimal generator. \par
Then, $t\mapsto p_{t}^{X}$ is a weak solution, in the sense of distributions, of the Kolmogorov forward equation:
\begin{equation}\label{4}
\left\{
\begin{array}{ll}
\frac{\partial p_{t}^{X}}{\partial t}=\mathcal{L}_{t}^{*}.p_{t}^{X} \\ \\
p_{0}^{X}=\delta_{0}
\end{array}
\right.
\end{equation}

\end{thm}

\begin{proof}

For any fixed $t\geq0$, we have $X_{t}\stackrel{(law)}{=}L_{t}$, i.e. $p_{t}^{X}=p_{t}^{L}.$ Now, applying the Fokker-Planck equation for Lévy processes,
$p_{t}^{L}=p_{t}^{X}$ is a weak solution of the Kolmogorov forward equation \eqref{4} described in the above theorem, so we can conclude.

\end{proof}

\begin{rem}

Considering an open set  $U=( 0,+\infty )\times\mathbb{R}$, the uniqueness of solution of the following problem in
$\mathcal{D'} (U)$,
\begin{equation*}
\left\{
\begin{array}{ll}
\frac{\partial u}{\partial t}=\mathcal{L}_{t}^{*} u  \\ \\
u_{0}=\delta_{0}
\end{array}
\right.
\end{equation*}
is an open question. In particular, we don't know if all the generalized solutions of this problem are densities of probability like $p_{t}^{X}$.
For the case without jumps, it's a result of M. Pierre, more details could be found in the monograph of F. Hirsch \textit{et al} \cite{Hirschbook}.  \par

\end{rem}

The previous theorem which is a classical result, can be meaningful if for a given IDT process, we know explicitly the associated Lévy process. 
Hence, the main purpose of this paper, is to exhibit numerous examples of IDT processes and their associated Lévy processes. 
Now, let us give the organization of the paper. \par
In Section 2, we give examples that motivated us to make construction of IDT processes and associated Lévy processes via Lévy sheet. 
In Section 3, we take up ideas of the monograph from F. Hirsch \textit{et al} \cite{Hirschbook} (Lévy, Sato and Gaussian sheets methods), 
and applying them to the framework of IDT processes, we present our sheets method. In Section 4, relying on the fact that IDT processes 
have the same one-dimensional marginals than Lévy processes, we give a weak Itô formula for IDT processes like the one for Gaussian processes 
given by F. Hirsch, B. Roynette and M. Yor \cite{Roynette}.\par
Following K. Es-Sebaiy and Y. Ouknine \cite{Ouknine}, we extend in Section 5 some results on Lévy processes presented in 
O. E. Barndorff-Nielsen, M. Maejima and K. Sato \cite{K.Sato}, to the case of IDT processes and we give a link between IDT processes and 
selfdecomposability. \par
In Section 6, we introduce a new concept of weak IDT process, which is basically asking for equation \eqref{1} to be satisfied only for 
one-dimensional marginals, i.e. an equality in law of random variables instead of an equality in law of processes. We also define the notion 
of 1-Lévy process, in the same spirit as 1-martingale in \textit{Peacocks and associated martingales}, and show that both notions are equivalent. \par
The last Sections 7 and 8 are devoted to the study of multiparameter IDT processes. Particularly, in Section 7 we invest multiparameter 
IDT processes introduced by K. Es-Sebaiy and Y. Ouknine \cite{Ouknine}, and in Section 8 we give a new approach of multiparameter IDT processes 
for which the $\mathbb{R}_{+}^{N}-$parameter Lévy process studied by O. E. Barndorff-Nielsen \textit{et al} \cite{Barndorff}, and the 
Lévy's $\mathbb{R}^{M}$-parameter Brownian motion studied by P. Lévy \cite{Levy}, N. N. Chentsov \cite{Chentsov} and H. P. Mckean Jr. \cite{Mckean}, 
are typical examples. To avoid confusions in the sequel, we are going to refer multiparameter IDT in the sense of 
K. Es-Sebaiy and Y. Ouknine \cite{Ouknine} as multiparameter IDT of type 1, and the multiparameter IDT in our sense as multiparameter IDT of type 2.

\section{The guiding example}

In this section, we use the notion of Lévy sheet (for which Brownian sheet is a special case) to construct Lévy processes associated to some given IDT 
processes. We refer to R. C. Dalang and J. B. Walsh \cite{Dalang} for definition of a Lévy sheet (and also Brownian sheet), and we recall the following result which would be used in the sequel.

\begin{thm}

Let $L=(L_{t}; t\geq 0)$ be an $\mathbb{R}^{d}$-valued Lévy process starting from 0.
Then, there exists an $\mathbb{R}^{d}$-valued two-parameter process $\tilde{L}=(\tilde{L}_{s,t}; s\geq0, t\geq0)$ satisfying the following properties:
\begin{itemize}
\item $\forall s,t \geq 0$, $\tilde{L}_{s,0}=\tilde{L}_{0,t}=0$.
\item Almost surely, for any $s,t\geq0$, $\tilde{L}_{s,.}$ and $\tilde{L}_{.,t}$ are càdlàg functions on $\mathbb{R}_{+}$.
\item Let, for $t\geq0$, $\mathcal{L}_{t}=\sigma(\tilde{L}_{u,v}; u\geq0, 0\leq v\leq t )$.
Then, for $0\leq t_{1} \leq t_{2}$, the process $(\tilde{L}_{s,t_{2}}-\tilde{L}_{s,t_{1}}; s\geq 0)$ is a Lévy process starting from 0, independent of
$\mathcal{L}_{t_{1}}$, which is distributed as $(L_{(t_{2}-t_{1})s}; s\geq 0 )$.
\item The two-parameter processes $(\tilde{L}_{s,t}; s,t\geq0)$ and $(\tilde{L}_{t,s}; s,t\geq0)$ have the same law.
\end{itemize}
The stochastic process $\tilde{L}$ is called the Lévy sheet extending the Lévy process L. Its law is fully determined by the one of L.

\end{thm}

\begin{proof}

A proof of this theorem may be found in R. C. Dalang and J. B. Walsh $\cite{Dalang}$, who, themselves, refer R. J. Adler \textit{et al} $\cite{Adler}$.

\end{proof}

We consider the standard Brownian motion $B=(B_{t};t \geq 0)$ and $\varphi$ a function of $L^{2}([0,1])$.
It's clear by R. Mansuy ($\cite{Mansuy}$ Example 3.4.) that the stochastic process $G^{\varphi}$ well-defined for all $t\geq0$ by $$G^{\varphi}_{t}=\int_{0}^{t} \varphi (\frac{u}{t})dB_{u}=\int_{0}^{1} \varphi (v)d_{v}B_{vt}$$
is an IDT Gaussian process. We try to construct a Lévy process having the same one-dimensional marginals distributions as $G^{\varphi}$. For this, let us consider $W=(W_{v,t}; v \geq 0, t \geq 0)$ the Brownian sheet extending 
the Brownian motion B and let $\bar{G}^{\varphi}$ the stochastic process well-defined for all $t\geq0$ by
$$ \bar{G}^{\varphi}_{t}=\int_{0}^{1} \varphi (v) d_{v}W_{v,t}.$$
Since $\forall t\geq 0 \mbox{ fixed }, \,\,\,(B_{ut}, u \geq 0)\stackrel{(law)}{=}(W_{u,t}, u \geq 0),$ the following result holds.

\begin{prop}

The stochastic process $\bar{G}^{\varphi}$ is an associated Gaussian Lévy process to the Gaussian IDT process $G^{\varphi}$, that is
$$ \forall t \geq 0 \mbox{ fixed },\,\, G^{\varphi}_{t}\stackrel{(law)}{=}\bar{G}^{\varphi}_{t}.$$

\end{prop}

In the sequel, we generalize the previous result to Lévy processes. First, let us give this remark due to K. Sato (\cite{Sato2004}, page 230).

\begin{rem}

We recall that $\int_{0}^{+\infty} f(s) dZ_{s}$ with $(Z_{t};t\geq0)$ a Lévy process, is defined as the limit in probability of
$\int_{0}^{h} f(s)dZ_{s}$ as $h\rightarrow+\infty$.

\end{rem}

\begin{prop}

Suppose that $(L_{t};t\geq0)$ is a Lévy process on $\mathbb{R}^{d}$, $f(s)$ a locally bounded function on $[0,+\infty)$ such that 
$\int_{0}^{+\infty} f(s)dL_{s}$ is well defined. \\
Then, the process $X=(X_{t};t\geq0)$ defined by:
$$X_{t}=\int_{0}^{+\infty} f(s)dL_{ts}$$
is an IDT process and its associated Lévy process $(\tilde{X}_{t}; t\geq0)$ is given by
$$\tilde{X}_{t}=\int_{0}^{+\infty} f(s)d_{s}\tilde{L}_{t,s}$$
where $(\tilde{L}_{t,s}; t,s\geq0)$ is the extending Lévy sheet of the Lévy process Z.

\end{prop}

\begin{proof}

For all $n,m\in\mathbb{N}^{\ast}$, we have:
$$\mathbb{E}exp\{i\sum_{k=1}^{m}<\theta_{k}, X_{n t_{k}}>\}=\mathbb{E} exp\{i\sum_{k=1}^{m} <\theta_{k}, \int_{0}^{+\infty} f(s) d L_{n s t_{k}}>\}$$
Thanks to the IDT property of the Lévy process L, we get
$$\mathbb{E}exp\{i\sum_{k=1}^{m}<\theta_{k}, X_{n t_{k}}>\}=\mathbb{E} exp\{i \sum_{k=1}^{m} 
<\theta_{k}, \sum_{j=1}^{n}\int_{0}^{+\infty}f(s)dL^{(j)}_{st_{k}}>$$
where $L^{(j)}$ are independent copies of L.\\
This implies that 
$$\mathbb{E}exp\{i\sum_{k=1}^{m}<\theta_{k}, X_{n t_{k}}>\}=(\mathbb{E} exp\{i \sum_{k=1}^{m} 
<\theta_{k}, \int_{0}^{+\infty}f(s)dL_{st_{k}}>)^{n}.$$
Hence, X is an IDT process. Now, according to the properties of the Lévy sheet $\tilde{L}$ given in Theorem 3.1, it is straightforward to prove that 
$\tilde{X}$ is a Lévy process and for any fixed $t\geq0$, we have $X_{t}\stackrel{(law)}{=}\tilde{X}_{t}$.

\end{proof}

\section{A general framework involving measurable sheet: Lévy sheet, Sato sheet, Gaussian sheet.}

In this section, we choose as in F. Hirsch \textit{et al} (\cite{Roynette} and \cite{HirschYor}), an adequate measurable sheet 
(Gaussian sheet, Sato sheet, Lévy sheet) from which we construct general IDT processes and associated Lévy processes.

\begin{thm}

Let $\Gamma$ be a measurable space and $\mu$ a "good measure" on $\Gamma$, i.e. a measure such that the integrals in \eqref{6} and \eqref{7} are well 
defined. For any $t\geq 0$, consider a real valued measurable process $(X_{\gamma,t};\gamma \in \Gamma)$.
Denote by D the usual Skorohod space of càdlàg functions and assume that the process $(X_{.,t};t\geq0)$ is a D-valued stochastically continuous IDT process.
Assume now the existence of a measurable sheet $(\tilde{X}_{\gamma ,t};\gamma \in \Gamma,t\geq 0)$ such that:\\ \\
(H1) For every $t\geq0$, $$X_{.,t}\stackrel{(law)}{=}\tilde{X}_{.,t}.$$
(H2) For all $0\leq s\leq t$,
$$\tilde{X}_{.,t}-\tilde{X}_{.,s} \mbox{ is independent of } \sigma(\tilde{X}_{\gamma,u};\gamma\in\Gamma, 0\leq u\leq s).$$
Then the process
\begin{equation}\label{6}
X^{ \mu }_{t}:=\int_{\Gamma} \mu (d\gamma)X_{\gamma,t},  \mbox{  for all  } t\geq 0,
\end{equation}
is an IDT process and its associated Lévy process is given by
\begin{equation}\label{7}
\tilde{X}^{ \mu }_{t}:=\int_{\Gamma} \mu (d\gamma )\tilde{X}_{ \gamma ,t},  \mbox{  for all } t\geq 0.
\end{equation}

\end{thm}

\begin{proof}

For every $m,n\geq 1$, $\theta =(\theta_{1},\cdots,\theta_{m})\in {\mathbb{R}}^{m}$, we have:
\begin{displaymath}
J(n, \theta):=\mathbb{E} exp \{ i\sum_{k=1}^{m} <\theta_{k}, X^{\mu}_{n t_{k}}> \}=\mathbb{E} exp \{ i\sum_{k=1}^{m} <\theta_{k}, \int_{\Gamma}
\mu (d \gamma )X_{ \gamma, n t_{k}}> \}.
\end{displaymath}
Using the fact that $(X_{.,t})_{t\geq0}$ is a D-valued IDT process, we have
\begin{displaymath}
J(n, \theta)=\mathbb{E} exp \{i \sum_{k=1}^{m} <\theta_{k}, \sum_{j=1}^{n} \int_{ \Gamma } \mu (d \gamma )X^{(j)}_{ \gamma , t_{k}}> \}.
\end{displaymath}
where $(X^{(j)}_{.,t})_{t\geq0}$ are independent copies of $(X_{.,t})_{t\geq0}$.\\
Then,
\begin{displaymath}
J(n, \theta)=\mathbb{E} \prod_{j=1}^{n} exp \{ i\sum_{k=1}^{m} <\theta_{k}, \sum_{j=1}^{n} \int_{ \Gamma} \mu (d \gamma )X^{(j)}_{ \gamma, t_{k}}> \}.
\end{displaymath}
According to the independence of the copies of $(X_{.,t})_{t\geq0}$, we have
\begin{displaymath}
J(n, \theta)= \prod_{j=1}^{n} \mathbb{E} exp \{ i\sum_{k=1}^{m} <\theta_{k}, \sum_{j=1}^{n} \int_{\Gamma} \mu (d\gamma )X^{(j)}_{ \gamma , t_{k}}> \}.
\end{displaymath}
Therefore we have
\begin{displaymath}
J(n, \theta):=(\mathbb{E} exp \{ i\sum_{k=1}^{m} <\theta_{k}, X^{\mu}_{t_{k}}> \})^{n}.
\end{displaymath}
Now, it is enough to prove that the D-valued process $(\tilde{X}_{.,t};t\geq0)$ has stationary increments. This follows easily from (H1), (H2) and the
stochastic continuity of the D-valued IDT process $(X_{.,t};t\geq0)$. And then $(\tilde{X}_{.,t};t\geq0)$ is a D-valued Lévy process.
Hence, it is straightforward to see that $(\tilde{X}^{ \mu }_{t};t\geq0$ is a Lévy process. \\
The proof is achieved.

\end{proof}

In the following, we are going to illustrate the above theorem by some examples which involve Lévy sheet, Sato sheet and Gaussian sheet.

\subsection{Lévy sheet}

We consider $\Gamma =\mathbb{R}_{+}$, a Lévy process $L=(L_{t}; t\geq0)$ and $\tilde{L}=(\tilde{L}_{\gamma, t}; \gamma,t\geq0)$ its extending Lévy sheet.
Now, setting $X_{\gamma, t}=L_{\gamma t}$ and $\tilde{X}_{\gamma, t}=\tilde{L}_{\gamma, t}$, the hypotheses (H1) and (H2) are satisfied.
We also consider a compactly supported measure $\mu$ on $\Gamma =\mathbb{R}_{+}$.
Then, the following result holds.

\begin{prop}

Assume that for any $t\geq0$ fixed,
$$\int_{0}^{+\infty}\mu(d\gamma)\mathbb{E}(|L_{\gamma t}|)=\int_{0}^{+\infty}\mu(d\gamma)\mathbb{E}(|\tilde{L}_{\gamma,t}|)<+\infty.$$
Then, the stochastic process $X^{\mu}$ defined by
$$X^{\mu}_{t}=\int_{0}^{+\infty} \mu (d\gamma)X_{\gamma,t}=\int_{0}^{+\infty} \mu (d\gamma)L_{\gamma t}$$
is an IDT process and its associated Lévy process is the process $\tilde{X}^{\mu}$ defined by
$$\tilde{X}^{\mu}_{t}=\int_{0}^{+\infty} \mu (d\gamma)\tilde{X}_{\gamma,t}=\int_{0}^{+\infty} \mu (d\gamma)\tilde{L}_{\gamma,t}.$$

\end{prop}

Now, in the sequel we calculate the Lévy measure and the Lévy exponent of the Lévy process $\tilde{X}^{\mu}$ via those of the Lévy sheet $\tilde{L}$.

\begin{rem}

We notice that, for all $\lambda\in\mathbb{R}$
\begin{displaymath}
\mathbb{E}e^{i\lambda \tilde{L}_{1,1}}=\mathbb{E}e^{i\lambda L_{1}}=e^{\psi (\lambda)},
\end{displaymath}
where $\psi$ is the Lévy symbol of the infinitely divisible random variable $L_{1}$. 
So, we call $\psi$ the Lévy symbol of the infinitely divisible random variable $\tilde{L}_{1,1}$.\\
According to the fact that for all $\lambda\in\mathbb{R}$ and for all $s,t\geq 0$, we have
\begin{displaymath}
\mathbb{E}e^{i\lambda \tilde{L}_{s,t}}=\mathbb{E}e^{i\lambda L_{st}}=e^{st \psi ( \lambda )}.
\end{displaymath}
We call the characteristic exponent of a given Lévy sheet $(\tilde{L}_{s,t}; {s,t \geq 0})$, the Lévy symbol of the infinitely divisible 
random variable $\tilde{L}_{1,1}.$

\end{rem}

\begin{prop}

Let $(\tilde{L}_{s,t}; s,t \geq 0)$ the precedent Lévy sheet and $\psi$ its Lévy exponent.
Denote $\psi^{(\mu)}$ the characteristic exponent of the Lévy process $(\tilde{X}_{t}^{\mu}, t\geq0)$.\\
Then,
\begin{equation}\label{8}
\psi^{(\mu)}(\lambda)=\int_{0}^{+\infty}\psi(\lambda\mu([h,+\infty))) dh.
\end{equation}

\end{prop}

\begin{proof}

First, let us remark that a Lévy process and its extending Lévy sheet have the same Lévy exponent. Then we have,
$$\mathbb{E}e^{i\lambda \tilde{L}_{s,t}}=\mathbb{E}e^{i\lambda L_{st}}=e^{st \psi( \lambda)}.$$
In other hand, we also have
$$\mathbb{E}e^{i\lambda \tilde{X}_{1}^{\mu}}=\mathbb{E}e^{i\lambda \int_{0}^{+\infty}\mu(ds)\tilde{X}_{s,1}}=
\mathbb{E}e^{i\lambda \int_{0}^{+\infty}\mu(ds)\tilde{L}_{s,1}}=\mathbb{E}e^{i\lambda\int_{0}^{+\infty}\mu(ds)L_{s}}.$$
Since $\mu$ is a compactly supported measure, we have $[\int_{h}^{+\infty} \mu (ds)L_{h}]_{0}^{+\infty}=0,$
and then, by integration by part we obtain $$\int_{0}^{+\infty}\mu (ds) L_{s}=\int_{0}^{+\infty}\mu([h,+\infty))dL_{h}.$$
Then,
$$\mathbb{E}e^{i\lambda \tilde{X}_{1}^{\mu}}=\mathbb{E} e^{i\lambda \int_{0}^{+\infty} \mu ([h,+\infty)) dL_{h}}$$
Now, let $H>0$ such that $Supp(\mu)\subset[0,H]$, then we get:
$$\mathbb{E}e^{i\lambda \tilde{X}_{1}^{\mu}}=\mathbb{E} e^{i\lambda \int_{0}^{H} \mu ([h,+\infty)) dL_{h}}.$$
For all $N\in\mathbb{N}^{\ast}$, let $(h_{j}=j\frac{H}{N})_{0\leq j\leq N}$ be a regular subdivision of $[0,H]$. \\
Then, 
$$\mathbb{E}e^{i\lambda \tilde{X}_{1}^{\mu}}=\mathbb{E}\lim_{N\rightarrow+\infty}e^{i\lambda\sum_{j=0}^{N-1}(L_{h_{j+1}}-L_{h_{j}})\mu([h_{j},+\infty))}$$
Thanks to the dominated convergence theorem, it follows that 
$$\mathbb{E}e^{i\lambda \tilde{X}_{1}^{\mu}}=\lim_{N\rightarrow +\infty}\mathbb{E}e^{i\lambda\sum_{j=0}^{N-1}(L_{h_{j+1}}-L_{h_{j}})\mu([h_{j},+\infty))}$$
And that is 
$$\mathbb{E}e^{i\lambda \tilde{X}_{1}^{\mu}}=\lim_{N\rightarrow +\infty}\mathbb{E}\prod_{j=0}^{N-1}e^{i\lambda (L_{h_{j+1}}-L_{h_{j}})
\mu ([h_{j},+\infty ))}$$
According to the fact that L is a Lévy process, we get
$$\mathbb{E}e^{i\lambda \tilde{X}_{1}^{\mu}}=\lim_{N\rightarrow +\infty}\prod_{j=0}^{N-1}
\mathbb{E}e^{i\lambda \mu ([ h_{j}, +\infty))(L_{h_{j+1}-h_{j}})}$$
This implies that 
$$\mathbb{E}e^{i\lambda \tilde{X}_{1}^{\mu}}=\lim_{N\rightarrow +\infty}\prod_{j=0}^{N-1}e^{(h_{j+1}-h_{j}) \psi (\lambda\mu([h_{j},+\infty)))}$$
Now, since $Supp(\mu)\subset[0,H]$ we get
$$\mathbb{E}e^{i\lambda \tilde{X}_{1}^{\mu}}=e^{\int_{0}^{H}\psi(\lambda \mu ([h,+\infty)))dh}=e^{\int_{0}^{+\infty}\psi(\lambda\mu([h,+\infty)))dh}$$
Hence, according to the the fact $\tilde{X}^{\mu}$ is a Lévy process, i.e.
$\mathbb{E}e^{i\lambda \tilde{X}_{1}^{\mu}}=e^{\psi^{(\mu)}(\lambda)},$
we immediately deduce \eqref{8}.

\end{proof}

\begin{prop}

Assume that the previous Lévy process $L=(L_{t}; t\geq0)$ is now a pure jump process and consider $(\tilde{L}_{s,t}; s,t \geq 0)$ its extending Lévy sheet.
Denote $\nu$ the Lévy measure of $\tilde{L}$ and set $\nu^{\mu}$ the Lévy measure of the Lévy process $\tilde{X}^{\mu}$.\\
Then for any non-negative Borel function $f$, we have
\begin{equation}\label{9}
\int {\nu}^{\mu}(dy) f(y)=\int_{0}^{ +\infty } dh \int \nu (dx) f(\mu ([h,+\infty))x).
\end{equation}

\end{prop}

\begin{proof}

First, notice that a Lévy process and its extending Lévy sheet have the same Lévy measure. 
Now, for any fixed $t\geq0$, we have $$\int_{0}^{+\infty} \mu (ds) \tilde{L}_{s,t}=\tilde{X}_{t}^{\mu}\stackrel{(law)}{=}X_{t}^{\mu}=\int_{0}^{+\infty}\mu(ds)L_{st}.$$
In particular we have $$\tilde{X}_{1}^{\mu}\stackrel{(law)}{=}\int_{0}^{+\infty} \mu (ds) L_{s}.$$
By integration by parts, we obtain
$$\int_{0}^{+\infty} \mu (ds) L_{s}=\int_{0}^{+\infty}\mu([h,+\infty))dL_{h}.$$
Hence,
\begin{equation*}
\mathbb{E}[exp(-\lambda \int_{0}^{+\infty}\mu(ds) L_{s})]=exp(-\int_{0}^{+\infty}dh \int \nu(dx)(1-e^{-\lambda\mu([h,+\infty))x}))
\end{equation*}
This implies that
\begin{equation*}
\mathbb{E}[exp(-\lambda \tilde{X}_{1}^{\mu}]=exp(-\int_{0}^{+\infty}dh \int \nu(dx)(1-e^{-\lambda\mu([h,+\infty))x}))
\end{equation*}
from which we immediately deduce \eqref{9}.

\end{proof}

\subsection{Sato sheet}

Let $\Gamma=\mathbb{R}_{+}$ and set $X_{\gamma, t}=tL_{\gamma}$ where $(L_{\gamma}; \gamma\geq0)$ is a Lévy process starting from 0 and
assume that $L_{1}$ is a strictly 1-stable random variable. Notice that any stable random variable is also selfdecomposable. Consider
$(\tilde{S}_{\gamma, t}; \gamma\geq0, t\geq0)$ the Sato sheet attached to $L_{1}$, i.e. the process which is characterized by:
\begin{itemize}
\item The process $(\tilde{S}_{.,t}; t\geq 0)$ is a D-valued process with independent increments.
\item For any fixed $t\geq0$, $(\tilde{S}_{.,t}; t\geq 0)$ is a D-valued 1-selfsimilar process i.e.
      $$\forall c>0 \,\,\, (\tilde{S}_{.,ct}; t\geq 0)\stackrel{(law)}{=}(c\tilde{S}_{.,t}; t\geq 0).$$
\item For any fixed $t\geq 0$, $\tilde{S}_{.,t}=(\tilde{S}_{\gamma, t}; \gamma\geq0)$ is a Lévy process.
\item $\tilde{S}_{1,1}\stackrel{(law)}{=}L_{1}$.
\end{itemize}
Setting $\tilde{X}_{\gamma, t}=\tilde{S}_{\gamma, t}$, it's clear that hypotheses (H1) and (H2) are satisfied.
Then, the following result holds.

\begin{prop}

Consider a compactly supported measure $\mu$ defined on $\mathbb{R}_{+}$, and assume that for any fixed $t\geq0$,
$$\int_{0}^{+\infty}t\mu(d\gamma)\mathbb{E}(|L_{\gamma}|)=\int_{0}^{+\infty}\mu(d\gamma)\mathbb{E}(|\tilde{S}_{\gamma,t}|)<+\infty.$$
Then
$$X^{\mu}_{t}=\int_{\mathbb{R}_{+}} \mu (d\gamma) X_{\gamma, t}=\int_{\mathbb{R}_{+}} \mu (d\gamma)tL_{\gamma}$$
is an IDT process and its associated Lévy process is
$$\tilde{X}^{\mu}_{t}=\int_{\mathbb{R}_{+}}\mu (d\gamma)\tilde{X}_{\gamma,t}=\int_{\mathbb{R}_{+}}\mu (d \gamma)\tilde{S}_{\gamma,t}.$$

\end{prop}

\subsection{Gaussian sheet}

We illustrate now our sheet method via Gaussian sheet. We refer to F. Hirsch \textit{et al} \cite{Hirschbook} for backgrounds on Gaussian sheet.

\begin{prop}

Let $(G_{\gamma, t}; \gamma\in\Gamma, t\geq0)$ be a family of centered Gaussian processes, that is for each $t\geq0$,
$G_{.,t}=(G_{\gamma, t}; \gamma\in\Gamma)$ is a centered Gaussian process. Consider $(\tilde{G}_{\gamma, t}; \gamma\in\Gamma, t\geq0)$ the Gaussian sheet
attached to this family of Gaussian processes, i.e. the process which satisfied in particular:
\begin{itemize}
 \item For any $t\geq0$ fixed, $(\tilde{G}_{\gamma, t}; \gamma\in\Gamma)\stackrel{(law)}{=}(G_{\gamma, t}; \gamma\in\Gamma)$.
 \item $(\tilde{G}_{., t}; t\geq0)$ is a D-valued process with independent increments.
\end{itemize}
Assume that $(G_{., t}; t\geq0)$ is a D-valued IDT process and let $\mu$ be a "good measure" i.e. a measure such that the following integrals
\eqref{ase} and \eqref{use}  are well-defined.
Then,
\begin{equation}\label{ase}
G^{\mu}_{t}:=\int_{\Gamma} \mu (d \gamma) G_{\gamma ,t}  \mbox{  for all } t\geq 0,
\end{equation}
is a Gaussian IDT process and its associated Gaussian Lévy process is given by
\begin{equation}\label{use}
\tilde{G}^{\mu}_{t}:=\int_{\Gamma} \mu (d \gamma) \tilde{G}_{\gamma ,t} \mbox{ for all } t\geq0.
\end{equation}

\end{prop}

\begin{exple}

We consider $\Gamma=\mathbb{R}_{+}$, $\mu$ a compactly supported measure on $\mathbb{R}_{+}$, and we set for any $\gamma,t \geq0$,
$$G_{\gamma, t}=\sqrt{t}B_{\gamma} \mbox{ and } \tilde{G}_{\gamma,t}=W_{\gamma, t}$$
where B and W are respectively the standard Brownian motion and the standard Brownian sheet. Then, the processes defined by
$$G^{\mu}_{t}:=\int_{\mathbb{R}_{+}}\mu (d\gamma)G_{\gamma ,t}=\int_{\mathbb{R}_{+}}\mu(d\gamma)\sqrt{t}B_{\lambda}\,\,\,  \mbox{  for all } t\geq0,$$
and
$$\tilde{G}^{\mu}_{t}:=\int_{\mathbb{R}_{+}}\mu(d\gamma)\tilde{G}_{\gamma ,t}=\int_{\mathbb{R}_{+}}\mu(d\gamma)W_{\gamma,t}\,\,\, \mbox{for all } t\geq0.$$
are respectively an IDT Gaussian process and its associated Gaussian Lévy process.

\end{exple}

\section{An integrated Itô formula for IDT processes}

As an introduction to this section, let us recall the following theorem due to F. Hirsch \textit{et al} ($\cite{Roynette}$ Theorem 2.1) for the simple case.

\begin{thm}

Let f be a $C^{2}$-function and $G=(G_{t}, t \geq 0)$ a centered Gaussian process.
Then,
\begin{equation*}
\mathbb{E}f(G_{t})=\mathbb{E}f(G_{0})+\int_{0}^{t}\mathbb{E}f'(G_{s})ds+\frac{1}{2}\int_{0}^{t}\mathbb{E}f"(G_{s})d_{s}Var(G_{s}).
\end{equation*}

\end{thm}

In what follows, this integrated Itô type formula, is extended to IDT processes.

\begin{thm}

Let $X=(X_{t};t\geq0)$ be a stochastically continuous IDT process and consider a $C^{2}$ real function $f$. We denote by $\nu$ the Lévy measure of the infinitely divisible random
variable $X_{1}$, and we also assume that $X_{1}$ is integrable. Then, we have the following weak Itô formula:
\begin{equation*}
\mathbb{E}f(X_{t})=f(0)+\mathbb{E}\int_{0}^{t} f'(X_{s-})d_{s}\mathbb{E}X_{s}+\frac{1}{2}\mathbb{E}\int_{0}^{t}f"(X_{s-})d_{s}Var(X_{s})
\end{equation*}
\begin{equation*}
+\mathbb{E}\int_{0}^{t}\int_{\mathbb{R}}(f(X_{s-}+x)-f(X_{s-})-x f'(X_{s-}))ds\nu(dx).
\end{equation*}

\end{thm}

\begin{proof}

Let $(L_{t};t\geq0)$ be the Lévy process which is associated to X. Then since $X_{1}\stackrel{(law)}{=}L_{1}$, $L_{1}$ is also integrable.
Denote the characteristic triplet of L by $(b,\sigma,\nu)$ and let $\mu$ be the random measure counting its jumps. 
Then, we have the following Lévy-Itô decomposition $$L_{t}=bt+\sigma W_{t}+\int_{0}^{t}\int_{\mathbb{R}}x(\mu(ds,dx)-ds\nu(dx)).$$
Now, applying Itô formula to the Lévy process L (see P. E. Protter \cite{Protter} Theorem 3.2), we get
\begin{equation*}
f(L_{t})=f(0)+\int_{0}^{t}f'(L_{s-})dL_{s}+\frac{1}{2}\int_{0}^{t}f"(L_{s-})d[L,L]^{c}_{s}
\end{equation*}
\begin{equation*}
 + \sum_{0\leq s\leq t} \{ f(L_{s})-f(L_{s-})-f'(L_{s-})\Delta (L_{s}) \}
\end{equation*}
That is 
\begin{equation*}
f(L_{t})=f(0)+\int_{0}^{t}f'(L_{s-})dL_{s}+\frac{1}{2}\int_{0}^{t}f"(L_{s-})d<L>_{s}
\end{equation*}
\begin{equation*}
+\int_{0}^{t}\int_{\mathbb{R}} \{ f(L_{s-}+x)-f(L_{s-})-xf'(L_{s-}) \}\mu(ds,dx)
\end{equation*}
Since, $\mathbb{E}f(X_{t})=\mathbb{E}f(L_{t})$ for any fixed $t\geq0$, it follows that
\begin{equation*}
\mathbb{E}f(X_{t})=f(0)+\mathbb{E}\int_{0}^{t}f'(L_{s-})dL_{s}+\frac{1}{2}\mathbb{E}\int_{0}^{t}f"(L_{s-})d<L>_{s}
\end{equation*}
\begin{equation*}
+\mathbb{E}\int_{0}^{t}\int_{\mathbb{R}}(f(L_{s-}+x)-f(L_{s-})-xf'(L_{s-}))ds\nu (dx).
\end{equation*}
Hence, we get
\begin{equation*}
\mathbb{E}f(X_{t})=f(0)+b\mathbb{E}\int_{0}^{t}f'(L_{s-})ds+\frac{\sigma^{2}}{2}\mathbb{E}\int_{0}^{t}f"(L_{s-})ds
\end{equation*}
\begin{equation*}
+\mathbb{E}\int_{0}^{t}\int_{\mathbb{R}}(f(L_{s-}+x)-f(L_{s-})-xf'(L_{s-}))ds\nu (dx).
\end{equation*}
Now, it follows that
\begin{equation*}
\mathbb{E}f(X_{t})=f(0)+\mathbb{E}\int_{0}^{t}f'(X_{s-})d_{s}\mathbb{E}X_{s}+\frac{1}{2}\mathbb{E}\int_{0}^{t}f"(X_{s-})
d_{s}Var(X_{s})
\end{equation*}
\begin{equation*}
+\mathbb{E}\int_{0}^{t}\int_{\mathbb{R}}(f(X_{s-}+x)-f(X_{s-})-xf'(X_{s-}))ds\nu (dx).
\end{equation*}

\end{proof}

\begin{rem}

The previous theorem also called weak Itô type formula for IDT processes, emphasize more again on the usefulness of constructing associated Lévy process 
for a given IDT process.

\end{rem}

\section{A link with selfdecomposability and related topics}

The notions of selfdecomposability, selfsimilarity, infinite divisibility and temporal selfdecomposability of processes, that we use in this section,
are in the sense of K. Sato \cite{Sato}. In order to give a link between IDT processes and selfdecomposability, let us first recall some results due 
to O. E. Barndorff-Nielsen, M. Maejima and K. Sato \cite{K.Sato}, which will be used in the sequel.

\begin{defn}

A stochastic process $X=(X_{t}; t\geq 0)$ on $\mathbb{R}^{d}$ is said to be of class $L_{1}$ or selfdecomposable (of order 1) if and only if, for every $c\in (0,1)$,
\begin{equation}\label{self}
X\stackrel{(law)}{=}c X' + U^{(c)}
\end{equation}
where $X'=(X'_{t}; t\geq 0)$ is a copy of $X$, $U^{(c)}=(U^{(c)}_{t}; t\geq 0)$ is an infinitely divisible process on $\mathbb{R}^{d}$, and X' and $U^{(c)}$ are independent.

\end{defn}

\begin{rem}

For all $m\in\mathbb{N}^{*}$, X is said to be of class $L_{m}$ or selfdecomposable of order m, if and only if, $U^{(c)}$ is of class $L_{m-1}$.

\end{rem}

\begin{thm}

If X is a selfdecomposable Lévy process (respectively a selfsimilar process) on $\mathbb{R}^{d}$, then, for every $c\in (0,1)$, the process $U^{(c)}$
in \eqref{self}, is also a Lévy process (respectively a selfsimilar process) on $\mathbb{R}^{d}$.

\end{thm}

\begin{proof}

The reader is referred to O. E. Barndorff-Nielsen, M. Maejima and K. Sato (\cite{K.Sato} Theorems 3.8 and 3.9).

\end{proof}

Analogous to the previous theorem, we have the following for IDT processes.

\begin{thm}

If X is a selfdecomposable IDT process on $\mathbb{R}^{d}$, then for any \\
$c\in (0,1)$, the process $U^{(c)}$ in \eqref{self} is an IDT process.

\end{thm}

\begin{proof}

Assume that X is a selfdecomposable IDT process. Then, we have
$$(X_{t}; t\geq 0)\stackrel{(law)}{=}(cX'_{t}+U^{(c)}_{t}; t\geq 0)$$
where X' is a copy of X and $U^{(c)}$ an infinitely divisible process independent of X'.
By the IDT property, we also have for all integer n,
$$ (X_{nt}; t\geq 0)\stackrel{(law)}{=}( X^{(1)}_{t} + \cdots + X^{(n)}_{t}; t\geq 0)$$
where $X^{(1)},\cdots,X^{(n)}$ are independent copies of X.\par
For all $m\in\mathbb{N}^{*}$, $\theta\in\mathbb{R}^{m}$, and by the selfdecomposability assumption, we have
\begin{equation*}
\mathbb{E} e^{i \sum_{k=1}^{m} <\theta_{k}, X_{t_{k}}>}=\mathbb{E} e^{i \sum_{k=1}^{m} <\theta_{k}, cX'_{t_{k}}+U^{(c)}_{t_{k}}>}.
\end{equation*}
According to the independence of X' and $U^{(c)}$, we have
\begin{equation*}
\mathbb{E} e^{i \sum_{k=1}^{m} <\theta_{k}, X_{t_{k}}>}=\mathbb{E} e^{i \sum_{k=1}^{m} <\theta_{k}, cX'_{t_{k}}>} \mathbb{E} e^{i \sum_{k=1}^{m}
<\theta_{k}, U^{(c)}_{t_{k}}>}.
\end{equation*}
Then,
\begin{equation*}
\mathbb{E} e^{i \sum_{k=1}^{m}<\theta_{k}, U^{(c)}_{t_{k}}>}=\frac{\mathbb{E} e^{i \sum_{k=1}^{m} <\theta_{k}, X_{t_{k}}>}}
{\mathbb{E} e^{i \sum_{k=1}^{m} <\theta_{k}, cX'_{t_{k}}>}}.
\end{equation*}
So, for all $n\in\mathbb{N}^{*}$,
\begin{equation*}
\mathbb{E} e^{i \sum_{k=1}^{m}<\theta_{k}, U^{(c)}_{nt_{k}}>}=\frac{\mathbb{E} e^{i \sum_{k=1}^{m} <\theta_{k}, X_{nt_{k}}>}}
{\mathbb{E} e^{i \sum_{k=1}^{m} <c\theta_{k}, X'_{nt_{k}}>}}.
\end{equation*}
Now, by the IDT property of X and X', we get
\begin{equation*}
\mathbb{E} e^{i \sum_{k=1}^{m}<\theta_{k}, U^{(c)}_{nt_{k}}>}=\frac{(\mathbb{E} e^{i \sum_{k=1}^{m} <\theta_{k}, X_{t_{k}}>})^{n}}
{(\mathbb{E} e^{i \sum_{k=1}^{m} <c\theta_{k}, X'_{t_{k}}>})^{n}}.
\end{equation*}
Thus,
\begin{equation*}
\mathbb{E} e^{i \sum_{k=1}^{m}<\theta_{k}, U^{(c)}_{nt_{k}}>}=(\mathbb{E} e^{i \sum_{k=1}^{m}<\theta_{k}, U^{(c)}_{t_{k}}>})^{n}.
\end{equation*}
Then, $U^{(c)}=(U^{(c)}_{t}; t\geq 0)$ is an IDT process.

\end{proof}

\begin{rem}

According to R. Mansuy (\cite{Mansuy} Section 4), any IDT process is an infinitely divisible process. The converse is not true. \\
\textbf{Counterexample}
\par
Let $V=(V_{t}; t\geq 0)$ be a stationary OU process on $\mathbb{R}^{d}$ defined by
$$V_{t}=\int_{0}^{t} e^{-\lambda (t-s)} dZ_{\lambda s}  \, \mbox{  for all  } t\geq 0,$$
where Z is a Lévy process with finite log-moments and $\lambda >0.$
Then, according to O. E. Barndorff-Nielsen, M. Maejima and K. Sato ($\cite{K.Sato}$ Theorem 4.1), V is an infinitely divisible process, but one can easily point that it is not an IDT process. \par
In fact, assume that V is an IDT process. Since V is stochastically continuous and according to R. Mansuy(\cite{Mansuy} Proposition 6.2), V is temporally selfdecomposable. However, according to O. E. Barndorff-Nielsen, M. Maejima and K. Sato (\cite{K.Sato} Theorem 5.12), V is not temporally selfdecomposable which contradict the fact that V is an IDT process.

\end{rem}

In the following, we give a link between IDT processes and selfdecomposability.

\begin{prop} \label{propG}

Let X be an $\mathbb{R}^{d}$-valued stochastically continuous IDT process and assume that X is a selfsimilar process. \\
Then, X is of class $L_{\infty}$ i.e. X is selfdecomposable of infinite order.

\end{prop}

\begin{proof}

Since X is a stochastically continuous IDT process, then according to R. Mansuy (\cite{Mansuy} Proposition 6.2),
X is temporally selfdecomposable (of infinite order). According to O. E. Barndorff-Nielsen, M. Maejima and K. Sato (\cite{K.Sato} Remark 5.9),
and the selfsimilarity assumption, we claim that X is selfdecomposable i.e. $$X\stackrel{(law)}{=}c X' + U^{(c)}$$ with X' and $U^{(c)}$ described
in \eqref{self}. Now thanks to Theorem 5.1 and Theorem 5.2, we show in the same way that $U^{(c)}$ is a selfdecomposable process
since it is a stochastically continuous and selfsimilar IDT process. Hence, one can conclude that X is selfdecomposable of infinite order,
either it is of class $L_{\infty}$.

\end{proof}

\begin{coro}

Let $L=(L_{t};t\geq0)$ be a Lévy process. Then,
\begin{equation*}
\mbox{ L is selfsimilar } \Rightarrow \mbox { L is selfdecomposable. }
\end{equation*}

\end{coro}

\begin{proof}

The corollary is direct consequence of the Proposition \ref{propG}.

\end{proof}

\section{Weak IDT processes}

This section is devoted to introduce the notion of weak IDT process which is deeply linked to the notion of 1-Lévy process.
We also give a sufficient condition for a weak IDT process to be a Lévy process and we prove that weak IDT processes could be obtained combining weak selfsimilarity and weak strict stability.

\begin{defn}

An $\mathbb{R}^{d}$-valued stochastic process $X=(X_{t};t\geq 0)$ is weakly IDT process if for all $n\in\mathbb{N}^{*}$, we have for any fixed $t\geq0$
\begin{equation}\label{5}
X_{nt}\stackrel{(law)}{=}X_{t}^{1} + \cdots + X_{t}^{n},
\end{equation}
where $X^{i}_{t}$, $i=1,\cdots,n$, are independent copies of $X_{t}$.

\end{defn}

\begin{exple}

All IDT processes are weak IDT processes. In particular, all Lévy processes are weak IDT processes.

\end{exple}

\begin{prop}

Let $0<\alpha\leq2$. A non-trivial, weakly strict $\alpha$-stable and weakly $(1/\alpha)$-selfsimilar process $(X_{t};t\geq0)$ is a weak IDT process.

\end{prop}

\begin{proof}

Firstly, X is weakly strict $\alpha$-stable, implies that for any fixed $t\geq0$, we have
$$n^{1/\alpha}X_{t}\stackrel{(law)}{=}X^{1}_{t}+\cdots +X^{n}_{t},$$
where $X^{1}_{t},\cdots, X^{n}_{t}$ are independent copies of $X_{t}$.\\
Secondly, X is weakly $(1/\alpha)$-selfsimilar, entails that for any fixed $t\geq0$, we have
$$X_{nt}\stackrel{(law)}{=}n^{1/\alpha}X_{t}.$$
Combining the two above equalities in law, the desired follows easily.

\end{proof}

\begin{prop}

Let $(X_{t};t\geq0)$ be a non-trivial, stochastically continuous, weak IDT process. Then $(X_{t};t\geq0)$ is weakly strict $\alpha$-stable if and
only if it is weakly $(1/\alpha)$-selfsimilar.

\end{prop}

\begin{proof}

The proof is straightforward, so we omit the details.

\end{proof}

\begin{prop}

Let $(X_{t}, t\geq0)$ be a stochastically continuous weak IDT process with independent increments. Then, $(X_{t}, t\geq0)$ is a Lévy process.

\end{prop}

\begin{proof}

Let $\theta\in\mathbb{R}$, $n,m\in\mathbb{N}^{\ast}$.
Since X is weakly IDT, we get by \eqref{5} that:
$$\mathbb{E} exp \{ i \theta X_{n}\}=(\mathbb{E} exp \{ i \theta X_{1}\})^{n}.$$
And
$$\mathbb{E} exp \{ i \theta X_{m/n}\}=(\mathbb{E} exp \{ i \theta X_{1}\})^{m/n}.$$
Then, by stochastic continuity of X and thanks to the density of $\mathbb{Q}_{+}$ in $\mathbb{R}_{+}$, we have for all $t>0$,
$$\mathbb{E} exp \{ i \theta X_{t}\}=(\mathbb{E} exp \{ i \theta X_{1}\})^{t}.$$
Now, for $0<s<t$, we have:
$$\mathbb{E} exp \{ i \theta X_{t-s}\}=(\mathbb{E} exp \{ i \theta X_{1}\})^{t-s}.$$
This implies that,
$$\mathbb{E} exp \{ i \theta X_{t-s}\}=\frac{\mathbb{E} exp \{ i \theta X_{t}\}}{\mathbb{E} exp \{ i \theta X_{s}\}}.$$
Thanks to the assumption of independent increments of X, we get:
$$\mathbb{E} exp \{ i \theta X_{t-s}\}=\mathbb{E} exp \{ i \theta (X_{t}-X_{s})\}.$$
From \eqref{5}, we have in particular
$$\mathbb{E} exp \{ i \theta X_{2.0}\}=(\mathbb{E} exp \{ i \theta X_{0}\})^{2}.$$
And then, $$\mathbb{E} exp \{ i \theta X_{0}\}=1.$$
Hence, it is easy to conclude that $X_{0}=0$ almost surely.
The proof is now complete.

\end{proof}

\begin{defn}

A stochastic process X is said to be a 1-Lévy process if there exist a Lévy process Y which is associated to X.

\end{defn}

So, we have the following which is just a reformulation of Proposition 1.2.

\begin{prop}

Any stochastically continuous IDT process is a 1-Lévy process.\\
The converse is not true.

\end{prop}

\textbf{Counterexample} \par
Consider the stochastic process X defined by
$$
X_{t}=
\left\{
\begin{array}{ll}
B_{t} & \mbox{  if  } t\leq \frac{1}{2} \\
B_{1/2}+(\sqrt{2}-1)B_{t-1/2} & \mbox{  if  } t > \frac{1}{2}
\end{array}
\right.
$$
\par
Then, according to H. Föllmer, C.T. Wu and M. Yor (\cite{Follmer} Remark 2.2), X is a continuous weak Brownian motion of order 1, so it is a 1-Lévy process.
In other hand, one can easily point that X is not an IDT (Gaussian) process. \par

\begin{prop}

Let $X=(X_{t};t\geq 0)$ be an $\mathbb{R}^{d}$-valued stochastic process which is assumed to be stochastically continuous.\\
Then, X is a weak IDT process if and only if it is a 1-Lévy process.

\end{prop}

\begin{proof}

The second implication is trivial. For the first, let us assume that X is a weak IDT process.
Then, we notice that $X_{1}$ is an infinitely divisible random variable and then there exist a Lévy process L (unique in law) such that
$X_{1}\stackrel{\textit{law}}{=}L_{1}$. Now, using the stochastic continuity of X and L, the result easily follows.

\end{proof}

\begin{coro}

Let $X=(X_{t}; t\geq0)$ be a 1-Lévy process which is provided to the independence of increments property. Assume that X is stochastically continuous. Then, X is a Lévy process (in law).

\end{coro}

\begin{proof}

The corollary is straightforward from Proposition 6.3 and Proposition 6.5. So the details are omitted.

\end{proof}

\begin{coro}

Let $X=(X_{t},t\geq 0)$ be a stochastically continuous weak IDT process.
Then, for all $s,t \geq0$, we have
$$X_{s+t}\stackrel{(law)}{=}X_{s}+X'_{t}$$
where $(X_{t}')$ is an independent copie of $(X_{t})$.

\end{coro}

\begin{proof}

The corollary is straightforward and a direct consequence of Proposition 6.5. So the details are omitted.

\end{proof}

\section{Multiparameter IDT processes of type 1}

In all the sequel, we denote $\delta(\textbf{a}):=\prod_{i=1}^{N} a_{i}$ for any $\textbf{a}=(a_{1},\cdots,a_{n})\in\mathbb{R}_{+}^{N}$, and we also 
consider a $(d \times d)$ invertible matrix Q (which for simplicity, is assumed to be symmetric), and for any $\alpha>0$ we let: 
$$\alpha^{Q}=e^{(\log{\alpha})Q}=\sum_{n=0}^{\infty} \frac{1}{n!} (\log{\alpha})^{n} Q^{n}.$$

\subsection{Definitions and some examples}

\begin{defn}

A multiparameter process $\{X(\textbf{s}); \textbf{s}=(s_{1},\cdots,s_{N})\in\mathbb{R}_{+}^{N}\}$ on $\mathbb{R}^{d}$ is a multiparameter IDT process of
type 1, if and only if, for all $\textbf{n}\in(\mathbb{N}^{\ast})^{N}$,
\begin{equation}
\{ X(\textbf{n.s}); \textbf{s} \in \mathbb{R}_{+}^{N} \} \stackrel{(law)}{=} \{ \sum_{i=1}^{\delta(\textbf{n})} X^{(i)}(\textbf{s}); \textbf{s} \in \mathbb{R}_{+}^{N} \},
\end{equation}
where $X^{(1)},\cdots,X^{(\delta(\textbf{n}))}$ are independent copies of X and $\textbf{n.s}:=(n_{1} s_{1},\cdots,n_{N} s_{N})$.

\end{defn}

\begin{exple}

(1) Let $\xi$ be a strictly $\alpha$-stable random variable, the process defined by
$$\{X(\textbf{s})=(s_{1}^{1/\alpha} \cdots s_{N}^{1/\alpha})\xi; \, \textbf{s}\in\mathbb{R}_{+}^{N}\}$$
is a multiparameter IDT process of type 1. \\ \\
(2) If X is a multiparameter IDT process of type 1 and $\mu$ a measure on $\mathbb{R}_{+}^{N}$ such that
$$X^{(\mu)}(\textbf{s})=\int_{\mathbb{R}_{+}^{N}} X(\textbf{u.s}) \mu (d\textbf{u}), \, \mbox{  } \textbf{s}\in\mathbb{R}_{+}^{N}$$
is well defined, then $X^{(\mu)}$ is also a multiparameter IDT of type 1. \\ \\
(3) Let $\{X(t); t\geq0\}$ be an IDT process, then the multiparameter process defined by
$$ Y(\textbf{s})=X(s_{1}s_{2} \cdots s_{N}) \mbox{ for any } \textbf{s}=(s_{1},\cdots,s_{N})\in\mathbb{R}_{+}^{N}$$
is a multiparameter IDT of type 1.

\end{exple}

In the following, we recall characterizations of Gaussian processes which are multiparameter IDT processes of type 1, due to K. Es-Sebaiy and Y. Ouknine \cite{Ouknine}.

\begin{prop}

Let $\{X(\textbf{s}); \textbf{s}\in\mathbb{R}_{+}^{N}\}$ be a stochastically continuous, centered Gaussian process. Then, the following properties are equivalent:

\begin{enumerate}
\item $\{X(\textbf{s}); \textbf{s}\in\mathbb{R}_{+}^{N}\}$ is a multiparameter IDT process of type 1.
\item The covariance function $\kappa(\textbf{s},\textbf{t}):=\mathbb{E}[X(\textbf{s})X(\textbf{t})]$, $\textbf{s}, \textbf{t} \in\mathbb{R}_{+}^{N}$,
satisfies
$$\kappa(\textbf{a.s},\textbf{a}.\textbf{t})=\delta(\textbf{a})\kappa(\textbf{s},\textbf{t}) \mbox{ for any } \textbf{a}\in(\mathbb{R}_{+}^{\ast})^{N}.$$
\item The process $\{X(\textbf{s}); \textbf{s}\in\mathbb{R}_{+}^{N}\}$ is $(N,d,\frac{1}{2})$-selfsimilar (in the sense of W. Ehm \cite{Ehm}), i.e.
$$\forall \textbf{a}\in(\mathbb{R}_{+}^{\ast})^{N}, \,\, \{X(\textbf{a.s}); \textbf{s}\in(\mathbb{R}_{+})^{N}\}\stackrel{(law)}{=}
\{\delta(\textbf{a})^{1/2}X(\textbf{s});\textbf{s}\in(\mathbb{R}_{+})^{N}\}.$$
\item Its Lamperti transform $\{Y(\textbf{y}):=e^{\{-\frac{1}{2}\sum_{j=1}^{N} y_{j}\}}X(e^{y_{1}},\cdots,e^{y_{N}}); \textbf{y}\in\mathbb{R}^{N}\}$ is a
strictly stationary process.
\end{enumerate}

\end{prop}

\begin{exple}

Let $\{B(\textbf{s}); \textbf{s}=(s_{1},\cdots,s_{N})\in\mathbb{R}_{+}^{N}\}$ be a Brownian sheet studied by S. Orey and W. E. Pruitt \cite{Orey} and many others, as the centered, real valued Gaussian random field with covariance function
$$\mathbb{E}[B(\textbf{s})B(\textbf{t})]=\prod_{i=1}^{N} s_{i} \wedge t_{i}.$$
It is easy to prove that $\{B(\textbf{s}); \textbf{s}=(s_{1},\cdots,s_{N})\in\mathbb{R}_{+}^{N}\}$ is a multiparameter Gaussian IDT process of type 1.

\end{exple}

\subsection{Links with operator stability and operator selfsimilarity}

In order to give a link between multiparameter IDT processes of type 1 and operator strict stability or operator selfsimilarity, let us recall the definitions of these concepts. These definitions are in the sense of W. Ehm \cite{Ehm}.

\begin{defn}

A multiparameter process $\{X(\textbf{s}); \textbf{s}\in\mathbb{R}_{+}^{N}\}$ is said $(N,d,Q)$-selfsimilar or Q-operator $(N,d)$-selfsimilar,
if and only if, for any $\textbf{a}\in(\mathbb{R}_{+}^{\ast})^{N}$,
\begin{equation}
\{X(\textbf{a}.\textbf{s}); \textbf{s}\in\mathbb{R}_{+}^{N}\}\stackrel{(law)}{=}\{\delta(\textbf{a})^{Q} X(\textbf{s}); \textbf{s}\in\mathbb{R}_{+}^{N}\}.
\end{equation}

\end{defn}

\begin{defn}

A multiparameter process $\{X(\textbf{s}); \textbf{s}\in\mathbb{R}_{+}^{N}\}$ is said strictly $(N,d,Q)$-stable or Q-operator strictly $(N,d)$-stable,
if and only if, for any $\textbf{n}\in(\mathbb{N}^{\ast})^{N}$
\begin{equation}
\{\sum_{i=1}^{\delta(\textbf{n})} X^{(i)}(\textbf{s}); \textbf{s}\in\mathbb{R}_{+}^{N}\}\stackrel{(law)}{=}\{\delta(\textbf{n})^{Q} X(\textbf{s}); \textbf{s}\in\mathbb{R}_{+}^{N}\},
\end{equation}
where $X^{(i)},\, i=1,\cdots,\delta(\textbf{n})$, are independent copies of X.

\end{defn}

Now, we are going to extend some results of K. Es-Sebaiy and Y. Ouknine \cite{Ouknine} on the connection between selfsimilarity and strict-stability
for one-parameter IDT to the case of multiparameter IDT of type 1.

\begin{prop}

A nontrivial, strictly $(N,d,Q)$-stable, $(N,d,Q)$-selfsimilar multiparameter process $\{X(\textbf{s}); \textbf{s}\in\mathbb{R}_{+}^{N}\}$ is a multiparameter IDT process of type 1.

\end{prop}

\begin{proof}

Firstly, since $\{X(\textbf{s}); \textbf{s}\in\mathbb{R}_{+}^{N}\}$ is strictly $(N,d,Q)$-stable, for all $\textbf{n}\in(\mathbb{N}^{\ast})^{N}$, we have
\begin{equation*}
\{\sum_{i=1}^{\delta(\textbf{n})} X^{(i)}(\textbf{s}); \textbf{s}\in\mathbb{R}_{+}^{N}\}\stackrel{(law)}{=}\{\delta(\textbf{n})^{Q} X(\textbf{s}); \textbf{s}\in\mathbb{R}_{+}^{N}\},
\end{equation*}
where $X^{(i)},\, i=1,\cdots,\delta(\textbf{n})$, are independent copies of X. \\
Secondly, since $\{X(\textbf{s}); \textbf{s}\in\mathbb{R}_{+}^{N}\}$ is $(N,d,Q)$-selfsimilar, for all $\textbf{n}\in(\mathbb{N}^{\ast})^{N}$, we have
\begin{equation*}
\{X(\textbf{n}.\textbf{s}); \textbf{s}\in\mathbb{R}_{+}^{N}\}\stackrel{(law)}{=}\{\delta(\textbf{n})^{Q} X(\textbf{s}); \textbf{s}\in\mathbb{R}_{+}^{N}\}.
\end{equation*}
Now, combining these two equalities (in law), it follows easily that X is a multiparameter IDT process of type 1.

\end{proof}

\begin{prop}

Let $\{X(\textbf{s}); \textbf{s}\in\mathbb{R}_{+}^{N}\}$ be a nontrivial, stochastically continuous, multi-parameter IDT process of type 1. Then, $\{X(\textbf{s}); \textbf{s}\in\mathbb{R}_{+}^{N}\}$ is strictly $(N,d,Q)$-stable if and only if it is $(N,d,Q)$-selfsimilar.

\end{prop}

\begin{proof}

$(\Longrightarrow)$ Since $\{X(\textbf{s}); \textbf{s}\in\mathbb{R}_{+}^{N}\}$ is a multi-parameter IDT process of type 1, for all $\textbf{n}=(n_{1},\cdots,n_{N})\in(\mathbb{N}^{\ast})^{N}$, we have
\begin{equation*}
\{ X(\textbf{n.s}); \textbf{s} \in \mathbb{R}_{+}^{N} \} \stackrel{(law)}{=} \{ \sum_{i=1}^{\delta(\textbf{n})} X^{(i)}(\textbf{s}); \textbf{s} \in \mathbb{R}_{+}^{N} \},
\end{equation*}
where $X^{(1)},\cdots,X^{(\delta(\textbf{n}))}$ are independent copies of X. \\
According to the assumption that X is strictly $(N,d,Q)$-stable, we get
\begin{equation*}
\{X(\textbf{n}.\textbf{s}); \textbf{s}\in\mathbb{R}_{+}^{N}\}\stackrel{(law)}{=}\{\delta(\textbf{n})^{Q} X(\textbf{s}); \textbf{s}\in\mathbb{R}_{+}^{N}\},
\end{equation*}
and also
\begin{equation*}
\{X(\textbf{q}.\textbf{s}); \textbf{s}\in\mathbb{R}_{+}^{N}\}\stackrel{(law)}{=}\{\delta(\textbf{q})^{Q} X(\textbf{s}); \textbf{s}\in\mathbb{R}_{+}^{N}\}
\end{equation*}
where $\textbf{q}=(q_{1},\cdots,q_{N})\in(\mathbb{Q}_{+}^{\ast})^{N}$ and $q_{j}=\frac{n_{j}}{m_{j}}$, $j=1,\cdots,N$, and $\textbf{n},\textbf{m}
\in(\mathbb{N}^{\ast})^{N}$. \\
Now, thanks to the stochastic continuity of X and the density of $(\mathbb{Q}_{+}^{\ast})^{N}$ in $(\mathbb{R}_{+}^{\ast})^{N}$, it follows that for all
$\textbf{a}=(a_{1},\cdots,a_{N})\in(\mathbb{R}_{+}^{\ast})^{N}$,
\begin{equation*}
\{X(\textbf{a}.\textbf{s}); \textbf{s}\in\mathbb{R}_{+}^{N}\}\stackrel{(law)}{=}\{\delta(\textbf{a})^{Q} X(\textbf{s}); \textbf{s}\in\mathbb{R}_{+}^{N}\}.
\end{equation*}

$(\Longleftarrow)$ Assume that $\{X(\textbf{s}); \textbf{s}\in\mathbb{R}_{+}^{N}\}$ is $(N,d,Q)$-selfsimilar, then we have for all $\textbf{a}=(a_{1},\cdots,a_{N})\in(\mathbb{R}_{+}^{\ast})^{N}$,
\begin{equation*}
\{X(\textbf{a}.\textbf{s}); \textbf{s}\in\mathbb{R}_{+}^{N}\}\stackrel{(law)}{=}\{\delta(\textbf{a})^{Q} X(\textbf{s}); \textbf{s}\in\mathbb{R}_{+}^{N}\}.
\end{equation*}
In particular, for all $\textbf{n}=(n_{1},\cdots,n_{N})\in(\mathbb{N}^{\ast})^{N}$, we have
\begin{equation*}
\{X( \textbf{n.s}); \textbf{s}\in\mathbb{R}_{+}^{N}\}\stackrel{(law)}{=}\{\delta(\textbf{n})^{Q} X(\textbf{s}); \textbf{s}\in\mathbb{R}_{+}^{N}\}.
\end{equation*}
Now, thanks to the fact that X is a multiparameter IDT process of type 1, it follows that
\begin{equation*}
\{\sum_{i=1}^{\delta(\textbf{n})} X^{(i)}(\textbf{s}); \textbf{s}\in\mathbb{R}_{+}^{N}\}\stackrel{(law)}{=}\{\delta(\textbf{n})^{Q} X(\textbf{s}); \textbf{s}\in\mathbb{R}_{+}^{N}\},
\end{equation*}
where $X^{(i)},\, i=1,\cdots,\delta(\textbf{n})$, are independent copies of X.\\
The proof is complete.

\end{proof}

\subsection{Links with operator selfdecomposability and temporal selfdecomposability}

In order to link multiparameter IDT processes of type 1 with operator selfdecomposability and temporal selfdecomposability, we give the definitions of these concepts.

\begin{defn}

A multiparameter process $X=\{X (\textbf{s}); \textbf{s}\in\mathbb{R}_{+}^{N}\}$ is said Q-operator (N,d)-selfdecomposable or (N,d,Q)-selfdecomposable,
if and only if, for any \\ $\textbf{c}=(c_{1},\cdots,c_{N})\in(0,1)^{N}$,
\begin{equation}\label{multiself}
X \stackrel{(law)}{=}\delta(\textbf{c})^{Q} X' + U^{(\textbf{c})},
\end{equation}
where $X'=\{X'(\textbf{s}); \textbf{s}\in\mathbb{R}_{+}^{N}\}$ is a copie of X and 
$U^{(\textbf{c})}=\{U^{(\textbf{c})}(\textbf{s});\textbf{s}\in\mathbb{R}_{+}^{N}\}$ a multiparameter infinitely divisible process independent of X'.

\end{defn}

\begin{defn}

An $\mathbb{R}^{d}$-valued multiparameter process $X=\{X(\textbf{s}); \textbf{s}\in\mathbb{R}_{+}^{N}\}$ is said (N,d)-temporally selfdecomposable, 
if and only if, for any $\textbf{c}\in (0,1)^{N}$,
there exist a process $U^{(\textbf{c})}=\{U^{(\textbf{c})}(\textbf{s}); \textbf{s}\in\mathbb{R}_{+}^{N}\}$ (called the $\textbf{c}$-residual of X) independent of $X$ and such that
\begin{equation}\label{multitempo}
\{X(\textbf{s}); \textbf{s}\in\mathbb{R}_{+}^{N}\}\stackrel{(law)}{=}\{X(\textbf{cs})+U^{(\textbf{c})}(\textbf{s}); \textbf{s}\in\mathbb{R}_{+}^{N}\}.
\end{equation}

\end{defn}

\begin{prop}

A stochastically continuous multiparameter IDT process of type 1 is (N,d)-temporally selfdecomposable of infinite order.

\end{prop}

\begin{proof}

The reader is referred to K. Es-Sebaiy and Y. Ouknine (\cite{Ouknine} Proposition 4.5) for a detailed proof.

\end{proof}

\begin{prop}

If $X=\{X(\textbf{s}); \textbf{s}\in\mathbb{R}_{+}^{N}\}$ is an (N,d,Q)-selfdecomposable multiparameter IDT process of type 1,
then for every $\textbf{c}\in(0,1)^{N}$, the process $U^{(\textbf{c})}$ defined in \eqref{multiself} is a multiparameter IDT process of type 1.

\end{prop}

\begin{proof}

Let $m\in\mathbb{N}^{\ast}$ and $(\theta_{1},\cdots,\theta_{m})\in\mathbb{R}^{m}$.
Since X is $(N,d,Q)$-selfdecomposable, for every $\textbf{c}\in(0,1)^{N}$, it follows by equality \eqref{multiself} that
$$\mathbb{E} exp\{i\sum_{k=1}^{m} <\theta_{k}, X(\textbf{s}^{k})>\}=\mathbb{E} exp\{i\sum_{k=1}^{m} <\theta_{k}, \delta(\textbf{c})^{Q} X'(\textbf{s}^{k})+U^{(\textbf{c})}(\textbf{s}^{k})>\}.$$
The independence of $U^{(\textbf{c})}$ and $X'$ implies that:
$$\mathbb{E} exp\{i\sum_{k=1}^{m} <\theta_{k}, X(\textbf{s}^{k})>\}=\mathbb{E} exp\{i\sum_{k=1}^{m} <\theta_{k}, \delta(\textbf{c})^{Q} X'(\textbf{s}^{k})>\} \mathbb{E} exp\{i\sum_{k=1}^{m} <\theta_{k}, U^{(\textbf{c})}(\textbf{s}^{k})>\}.$$
And then we get,
$$\mathbb{E} exp\{i\sum_{k=1}^{m} <\theta_{k}, U^{(\textbf{c})}(\textbf{s}^{k})>\}=\frac{\mathbb{E} exp\{i\sum_{k=1}^{m} <\theta_{k}, X(\textbf{s}^{k})>\}}{\mathbb{E} exp\{i\sum_{k=1}^{m} <\theta_{k}, \delta(\textbf{c})^{Q} X'(\textbf{s}^{k})>\}}.$$
Now, since X and X' are multiparameter IDT processes of type 1, this implies that for all $\textbf{n}=(n_{1},\cdots,n_{N})\in(\mathbb{N}^{\ast})^{N}$, we get
$$\mathbb{E} exp\{i\sum_{k=1}^{m} <\theta_{k}, U^{(\textbf{c})}(\textbf{n}.\textbf{s}^{k})>\}=\frac{(\mathbb{E} 
exp\{i\sum_{k=1}^{m} <\theta_{k}, X(\textbf{s}^{k})>\})^{\delta(\textbf{n})}}{(\mathbb{E} 
exp\{i\sum_{k=1}^{m} <\theta_{k}, \delta(\textbf{c})^{Q}X'(\textbf{s}^{k})>\})^{\delta(\textbf{n})}}.$$
Then, $U^{(\textbf{c})}$ is a multiparameter IDT process of type 1, that is:
$$\mathbb{E} exp\{i\sum_{k=1}^{m} <\theta_{k}, U^{(\textbf{c})}(\textbf{n}.\textbf{s}^{k})>\}
=(\mathbb{E} exp\{i\sum_{k=1}^{m} <\theta_{k}, U^{(\textbf{c})}(\textbf{s}^{k})>\})^{\delta(\textbf{n})}.$$

\end{proof}

\begin{coro}

A stochastically continuous, (N,d,Q)-selfsimilar multiparameter IDT process of type 1 is (N,d,Q)-selfdecomposable of infinite order.

\end{coro}

\begin{proof}

The corollary is a consequence of Proposition 7.2 and Proposition 7.3.

\end{proof}

\subsection{Study of the special case N=2}

In order to study the bi-parameter IDT processes of type 1, we recall some results on these processes.

\begin{prop}

Let $\{Z(\textbf{s}); \textbf{s}=(s_{1},s_{2})\in\mathbb{R}_{+}^{2}\}$ be a L\'{e}vy sheet on $\mathbb{R}^{d}$, then it is a bi-parameter IDT process of type 1.

\end{prop}

\begin{proof}

The reader is referred to K. Es-Sebaiy and Y. Ouknine (\cite{Ouknine} Proposition 4.1).

\end{proof}

\begin{prop}

Let $\{X(\textbf{s}); \textbf{s}\in\mathbb{R}_{+}^{2}\}$ be a stochastically continuous bi-parameter IDT process of type 1.
Then there exist a L\'{e}vy sheet $\{Z(\textbf{s}); \textbf{s}\in\mathbb{R}_{+}^{2}\}$ such that
$$X(\textbf{s})\stackrel{(law)}{=}Z(\textbf{s}) \mbox{ for any fixed } \textbf{s}=(s_{1},s_{2})\in\mathbb{R}_{+}^{2}.$$

\end{prop}

\begin{proof}

The reader is referred to K. Es-Sebaiy and Y. Ouknine (\cite{Ouknine} Proposition 4.3).

\end{proof}

In the one-parameter case, K. Es-Sebaiy and Y. Ouknine \cite{Ouknine} showed that under hypotheses of stochastic continuity and independence of increments, an IDT process is a L\'{e}vy process.
Also, they proved the inheritance of IDT property under time change when base processes are L\'{e}vy processes. \par
In the following, we are going to point out that those great properties are not shared for the L\'{e}vy sheet and the bi-parameter IDT processes of type 1.

\begin{exple}

Let $\{X(\textbf{u}); \textbf{u}=(u_{1},u_{2})\in\mathbb{R}_{+}^{2}\}$ be a L\'{e}vy sheet on $\mathbb{R}^{d}$ and denote by $\mu$ the law of $X(1,1)$. 
Let  $\{Z(\textbf{s}); \textbf{s}=(s_{1},s_{2})\in\mathbb{R}_{+}^{2}\}$ be a bi-parameter IDT process of type 1 such that $Z(\textbf{s})\in\mathbb{R}_{+}^{2}$ a.s. for all $\textbf{s}=(s_{1},s_{2})\in\mathbb{R}_{+}^{2}$.
Assume that these two processes are independent and define $Y(\textbf{s})=X_{Z(\textbf{s})}$.
We give below an example where $\{Y(\textbf{s}); \textbf{s}=(s_{1},s_{2})\in\mathbb{R}_{+}^{2}\}$ is not a bi-parameter IDT of type 1. \par
For this, we consider $Z(\textbf{s})=(s_{1}s_{2},s_{1}s_{2})$ for $\textbf{s}=(s_{1},s_{2})\in\mathbb{R}_{+}^{2}$, then it is obvious that
$\{Z(\textbf{s}); \textbf{s}=(s_{1},s_{2})\in\mathbb{R}_{+}^{2}\}$ is a bi-parameter IDT of type 1 on $\mathbb{R}_{+}^{2}$.
Then, we have $Y(s_{1},s_{2})=X(s_{1}s_{2},s_{1}s_{2})$ and $Y(s_{1},s_{2})$ has characteristic function $\hat{\mu}(z)^{(s_{1}s_{2})^{2}}$
for $(s_{1},s_{2})\in\mathbb{R}_{+}^{2}$.
If Y is a bi-parameter IDT of type 1, there may exist a L\'{e}vy sheet $\{L(\textbf{s}); \textbf{s}=(s_{1},s_{2})\in\mathbb{R}_{+}^{2}\}$ such that
$$Y(s_{1},s_{2})\stackrel{(law)}{=}L(s_{1},s_{2}) \mbox{ for any fixed } (s_{1},s_{2})\in\mathbb{R}_{+}^{2}.$$
And this implies that the characteristic function of $L(s_{1},s_{2})$ may be $\hat{\mu}(z)^{(s_{1}s_{2})^{2}}$ for $(s_{1},s_{2})\in\mathbb{R}_{+}^{2}$
where $\hat{\mu}$ is the Fourier transform of $\mu=\mathcal{L}(Y(s_{1},s_{2}))=\mathcal{L}(L(s_{1},s_{2}))$. Unfortunately, this is incompatible with the characteristic 
function of a L\'{e}vy sheet.

\end{exple}

\begin{rem}

Following J. Pedersen and K. Sato \cite{PedersenSato}, one may prove that $Y(s_{1},s_{2})$ need not even be infinitely divisible.
In fact, if we let $X(u_{1},u_{2})=u_{1}u_{2}$, which is a L\'{e}vy sheet on $\mathbb{R}_{+}$, then $Y(s_{1},s_{2})=Z^{1}(s_{1},s_{2})Z^{2}(s_{1},s_{2})$,
where $Z^{i}(s_{1},s_{2})$ denotes the ith coordinate of $Z(s_{1},s_{2})$ for $i=1,2$, and those coordinates are bi-parameter IDT of type 1. 
Since $Z^{i}(s_{1},s_{2})$, $i=1,2$, is an infinitely divisible random variable, we just have to construct $Z(s_{1},s_{2})$ such that the product of 
the coordinates is not infinitely divisible. Following D. N. Shanbhag \textit{et al} \cite{Shanbhag}, V. Rohatgi \textit{et al} \cite{Rohatgi}, such a 
construction is possible.

\end{rem}

\begin{exple}

Let $\{X(s_{1},s_{2}); (s_{1},s_{2})\in\mathbb{R}_{+}^{2}\}$ be a stochastically continuous bi-parameter IDT of type 1.
For any rectangle $B=(s_{1},t_{1}]\times(s_{2},t_{2}]$ in $\mathbb{R}_{+}^{2}$ such that $s_{i} \leq t_{i}$, we set
$$X(B)=X(t_{1},t_{2})-X(s_{1},t_{2})-X(t_{1},s_{2})+X(s_{1},s_{2}).$$
If $B_{1},\cdots,B_{n}$ are disjoint rectangles in $\mathbb{R}_{+}^{2}$ and $B=\cup_{j=1}^{n}B_{j}$, we set $$X(B)=\sum_{j=1}^{n} X(B_{j}).$$
For any $n\geq2$ and $B_{1},\cdots,B_{n}$ disjoint rectangles in $\mathbb{R}_{+}^{2}$, assume that $X(B_{1}),\cdots,X(B_{n})$ are independent.
For all $(s_{1},s_{2}), (t_{1},t_{2})\in\mathbb{R}_{+}^{2}$ such that $s_{i} \leq t_{i}$, we consider two rectangles $\tilde{B}=(0,t_{1}-s_{1}]\times(0,t_{2}-s_{2}]$ and $B=(s_{1},t_{1}]\times(s_{2},t_{2}]$. Then, we have $B=\tilde{B}+\textbf{s}$ for $\textbf{s}=(s_{1},s_{2})\in\mathbb{R}_{+}^{2}$. \\
Now, we shall give a stochastically continuous biparameter IDT process X such that $\mathcal{L}(X(B)) \neq \mathcal{L}(X(\tilde{B}))$.
For this, let $\{X(s_{1},s_{2}); (s_{1},s_{2})\in\mathbb{R}_{+}^{2}\}$ defined by
$$X(s_{1},s_{2})=s_{1}^{1/\alpha}s_{2}^{1/\alpha}\xi \mbox{ for all } (s_{1},s_{2})\in\mathbb{R}_{+}^{2},$$
where $\xi$ is a strictly $\alpha$-stable random variable of index $0<\alpha<2$.
We choose $(t_{1},t_{2})$ such that $t_{1}=2s_{1}$ and $t_{2}=2s_{2}$. Then, we get:
$$X(B)=X((s_{1},2s_{1}]\times(s_{2},2s_{2}])=(4^{1/\alpha}+1-2\times 2^{1/\alpha})s_{1}^{1/\alpha}s_{2}^{1/\alpha}\xi$$
and
$$X(\tilde{B})=X((0,s_{1}]\times(0,s_{2}])=s_{1}^{1/\alpha}s_{2}^{1/\alpha}\xi.$$
Denote by $\mu$ the law of $\xi$ and $\hat{\mu}$ its Fourier transform. Then we have
$$\mathbb{E}exp\{izX(B)\}=\mathbb{E}e^{iz(4^{1/\alpha}+1-2\times 2^{1/\alpha})s_{1}^{1/\alpha}s_{2}^{1/\alpha}\xi}
=[\hat{\mu}\{(4^{1/\alpha}+1-2\times 2^{1/\alpha})z\}]^{s_{1}s_{2}}$$
and
$$\mathbb{E}exp\{izX(\tilde{B})\}=\mathbb{E}e^{izs_{1}^{1/\alpha}s_{2}^{1/\alpha}\xi}
=[\hat{\mu}(z)]^{s_{1}s_{2}}.$$
Now, $\mathcal{L}(X(B))=\mathcal{L}(X(\tilde{B}))$ if and only if for all $z\in\mathbb{R}$, we have
$$[\hat{\mu}(z)]^{s_{1}s_{2}}=[\hat{\mu}\{(4^{1/\alpha}+1-2\times 2^{1/\alpha})z\}]^{s_{1}s_{2}}.$$
This implies that $4^{1/\alpha}-2\times 2^{1/\alpha}=0$ i.e. $\alpha=1$.
Hence, for $\alpha\neq1$, we get $$\mathcal{L}(X(B))\neq\mathcal{L}(X(\tilde{B})),$$
where B and $\tilde{B}$ are rectangles in $\mathbb{R}_{+}^{2}$ and $B=\tilde{B}+\textbf{s}$ with $\textbf{s}=(s_{1},s_{2})\in\mathbb{R}_{+}^{2}$.
And then, X is not a L\'{e}vy sheet.

\end{exple}

\section{Multiparameter IDT processes of type 2}

\subsection{Definition and motivating examples}

Following O. E. Barndorff-Nielsen \textit{et al} \cite{Barndorff}, we are going to construct a multiparameter process issued from an IDT process with 
independent components which naturally need to be IDT. But it will not be the case in the sense of K. Es-Sebaiy and Y. Ouknine \cite{Ouknine}. 
Fortunately, our new approach solve this interesting case.\\
Now, we consider N independent IDT processes $\{X_{1}(t)\},\cdots,\{X_{N}(t)\}$ respectively on $\mathbb{R}^{d_{1}},\cdots,\mathbb{R}^{d_{N}}$.
The stacked process $\{X(t)\}$ defined by $X(t)=(X_{1}(t),\cdots,X_{N}(t))^{\top}$ is then an IDT process on $\mathbb{R}^{d}$, where $d=d_{1}+\cdots+d_{N}$.
We deal with a multiparameter time $\textbf{s}=(s_{1},\cdots,s_{N})\in\mathbb{R}_{+}^{N}$ to define a multiparameter 
process $\{X(\textbf{s});\textbf{s}\in\mathbb{R}_{+}^{N}\}$ by 
$X(\textbf{s})=(X_{1}(s_{1}),\cdots,X_{N}(s_{N}))^{\top}$. 
Then, the process $\{X(\textbf{s}); \textbf{s}\in\mathbb{R}_{+}^{N}\}$ defined before, is an IDT in the following sense.

\begin{defn}

A multiparameter process $\{X(\textbf{s}); \textbf{s}=(s_{1},\cdots,s_{N})\in\mathbb{R}_{+}^{N}\}$ on $\mathbb{R}^{d}$ is said to be a multiparameter IDT process of type 2, if and only if, for all $n\in\mathbb{N}^{\ast}$,
\begin{equation}\label{16}
\{X(n.\textbf{s}); \textbf{s}\in\mathbb{R}_{+}^{N}\}\stackrel{(law)}{=}\{\sum_{i=1}^{n} X^{(i)}(\textbf{s}); \textbf{s}\in\mathbb{R}_{+}^{N}\},
\end{equation}
where $X^{(1)},\cdots,X^{(n)}$ are independent copies of X and $n.\textbf{s}:=(ns_{1},\cdots,ns_{N})$.

\end{defn}

\begin{exple}

(1) Let $\xi$ be a strictly $\alpha$-stable random variable with $0<\alpha<2$. Then, the process defined by
$$\{X(\textbf{s})=(s_{1}^{1/\alpha}+ \cdots +s_{N}^{1/\alpha})\xi; \, \textbf{s}\in\mathbb{R}_{+}^{N}\}$$
is a multiparameter IDT process of type 2. \\
\\
(2) If X is a multiparameter IDT process of type 2 and $\mu$ a measure on $\mathbb{R}_{+}^{N}$ such that
$$X^{(\mu)}(\textbf{s})=\int_{\mathbb{R}_{+}^{N}} X(\textbf{u.s}) \mu(d\textbf{u}); \, \mbox{  } \textbf{s}\in\mathbb{R}_{+}^{N}$$
is well defined, then $X^{(\mu)}$ is a multiparameter IDT process of type 2. \\
\\
(3) Let $\{Z(t); t\geq0\}$ be an IDT process on $\mathbb{R}^{d}$. Fix $\textbf{c}=(c_{1},\cdots,c_{N})\in\mathbb{R}_{+}^{N}$ and for any $\textbf{s}=(s_{1},\cdots,s_{N})\in\mathbb{R}_{+}^{N}$ define $$X(\textbf{s})=Z(<\textbf{c},\textbf{s}>)=Z(c_{1}s_{1}+\cdots+c_{N}s_{N}).$$
Then $\{X(\textbf{s});\, \textbf{s}\in\mathbb{R}_{+}^{N}\}$ is a multiparameter IDT process of type 2.

\end{exple}

In the following, we are going to characterize Gaussian processes which are multiparameter IDT of type 2. A kind example of these processes, is the L\'{e}vy's $\mathbb{R}^{M}$-parameter Brownian motion.

\begin{prop}

Let $\{G(\textbf{s}); \textbf{s}\in\mathbb{R}_{+}^{N}\}$ be a centered Gaussian multiparameter process, which is assumed to be continuous in probability. Then the following properties are equivalent:
\begin{enumerate}
\item $\{G(\textbf{s}); \textbf{s}\in\mathbb{R}_{+}^{N}\}$ is a multiparameter IDT process of type 2.
\item The covariance function $c(\textbf{s}^{1},\textbf{s}^{2}):=\mathbb{E}[G(\textbf{s}^{1})G(\textbf{s}^{2})]$, $0\preceq\textbf{s}^{1}\preceq\textbf{s}^{2}$, satisfies:
    $$\forall \alpha >0, \mbox{    } c(\alpha\textbf{s}^{1},\alpha\textbf{s}^{2})=\alpha c(\textbf{s}^{1},\textbf{s}^{2}).$$
\item The process $\{G(\textbf{s}); \textbf{s}\in\mathbb{R}_{+}^{N}\}$ is $\frac{1}{2}$-selfsimilar, i.e.
$$\forall a>0, \,\, \{X(a.\textbf{s}); \textbf{s}\in\mathbb{R}_{+}^{N}\}\stackrel{(law)}{=}\{a^{1/2}X(\textbf{s});\textbf{s}\in\mathbb{R}_{+}^{N}\}.$$
\end{enumerate}

\end{prop}

\begin{proof}

$(1\Leftrightarrow2)$ $\{G(\textbf{s});\textbf{s}\in\mathbb{R}_{+}^{N}\}$ is a multiparameter IDT process of type 2, if and only if, for any
$n\in\mathbb{N}^{\ast}$, $\textbf{s},\textbf{t}\in\mathbb{R}_{+}^{N}$,
$$c(n.\textbf{s}, n.\textbf{t})=n c(\textbf{s}, \textbf{t})$$
and also, for any $q\in\mathbb{Q}_{+}$, $\textbf{s}, \textbf{t}\in\mathbb{R}_{+}^{N}$,
$$c(q\textbf{s}, q\textbf{t})=q c(\textbf{s}, \textbf{t}).$$
Now by the stochastic continuity of G (i.e. the continuity of c), and the density of $\mathbb{Q}_{+}$ in $\mathbb{R}_{+}$, the desired result follows. \\
\\
$(2\Leftrightarrow3)$ Since the law of a centered Gaussian process is determined by its covariance function, we easily get that G is
$\frac{1}{2}$-selfsimilar.

\end{proof}

\begin{exple}

The L\'{e}vy's $\mathbb{R}^{M}$-parameter Brownian motion $\{X(\textbf{s});\textbf{s}\in\mathbb{R}^{M}\}$ is a Gaussian process characterized in law by
the properties that
$$\mathbb{E}[X(\textbf{s})]=0 \mbox{ and } \mathbb{E}[X(\textbf{s})X(\textbf{t})]=
\frac{1}{2}(||\textbf{s}||+||\textbf{t}||-||\textbf{t}-\textbf{s}||) \mbox{ for } \textbf{t,\,s}\in\mathbb{R}^{M}.$$
Hence $\mathcal{L}(X(\textbf{s}))=\mathcal{N}(0,||\textbf{s}||)$ and $\mathcal{L}(X(\textbf{t})-X(\textbf{s}))=\mathcal{N}(0,||\textbf{t-s}||)$.
Thus, the restriction $\{X(\textbf{s});\textbf{s}\in\mathbb{R}_{+}^{M}\}$ satisfies the previous proposition, then it is a multiparameter Gaussian IDT of
type 2.

\end{exple}

\subsection{Links with \texorpdfstring{$\mathbb{R}_{+}^{N}-$}{}parameter L\'{e}vy processes}

In order to point out links between $\mathbb{R}_{+}^{N}-$parameter IDT processes of type 2 and $\mathbb{R}_{+}^{N}-$parameter L\'{e}vy processes, 
we refer to J. Pederson and K. Sato \cite{PedersenSato} for definition and properties of $\mathbb{R}_{+}^{N}-$parameter 
L\'{e}vy processes.

\begin{prop}\label{prop8.1}

Let $X=\{X(\textbf{s});\textbf{s}=(s_{1},\cdots,s_{N})\in\mathbb{R}_{+}^{N}\}$ be a multiparameter IDT process of type 2.
For all $\textbf{s}=(s_{1},\cdots,s_{N})\in\mathbb{R}_{+}^{N}$, denote $\mu(\textbf{s})=\mathcal{L}(X(\textbf{s}))$ the law of the variable
$X(\textbf{s})$. Then, $\{\mu(\textbf{s}); \textbf{s}\in\mathbb{R}_{+}^{N}\}$ is a $\mathbb{R}_{+}^{N}-$parameter convolution semigroup on $\mathbb{R}^{d}$.

\end{prop}

\begin{proof}

From (\ref{16}), it is clear that for any $\textbf{s}\in\mathbb{R}_{+}^{N}$, $X(\textbf{s})$ is an infinitely divisible random variable i.e.
$\mu(\textbf{s})$ is an infinitely divisible distribution on $\mathbb{R}^{d}$. Now, for $j=1,\cdots,N$, let us consider
$e^{j}=(\delta_{jk})_{1\leq k\leq N}$,
where $\delta_{jk}=0$ or 1 according as $k \neq j$ or $k=j$. Set $\rho_{j}=\mu(e^{j})$ for $j=1,\cdots,N$.
Since $e^{j}\in\mathbb{R}_{+}^{N}$, $j=1,\cdots,N$, the distributions $\rho_{1},\cdots,\rho_{N}$ are infinitely divisible with characteristic triplets
$(A_{1},\nu_{1},\gamma_{1}),\cdots,(A_{N},\nu_{N},\gamma_{N})$ respectively. Then, following J. Pedersen and K. Sato (\cite{PedersenSato} Theorem 1.1),
$\{\mu(\textbf{s}); \textbf{s}\in\mathbb{R}_{+}^{N}\}$ is a $\mathbb{R}_{+}^{N}-$parameter convolution semigroup on $\mathbb{R}^{d}$. Also,
for any $\textbf{s}=s_{1}e^{1}+ \cdots + s_{N}e^{N}\in\mathbb{R}_{+}^{N}$, $\mu_{\textbf{s}}$ is infinitely divisible with characteristic triplet
$(A_{\textbf{s}}, \nu_{\textbf{s}}, \gamma_{\textbf{s}})$ such that
$$A_{\textbf{s}}=s_{1}A_{e^{1}}+ \cdots + s_{N}A_{e^{N}},$$ 
$$\nu_{\textbf{s}}=s_{1}\nu_{e^{1}}+ \cdots + s_{N}\nu_{e^{N}},$$
$$\gamma_{\textbf{s}}=s_{1}\gamma_{e^{1}}+ \cdots + s_{N}\gamma_{e^{N}}.$$

\end{proof}

\begin{thm}

Any $\mathbb{R}_{+}^{N}-$parameter L\'{e}vy process is a multiparameter IDT process of type 2. Conversely, any stochastically continuous multiparameter IDT process of type 2 with independent increments, is an $\mathbb{R}_{+}^{N}-$parameter L\'{e}vy process.

\end{thm}

\begin{proof}

$(\Longrightarrow)$ Let $\{X(\textbf{s});\textbf{s}=(s_{1},\cdots,s_{N})\in\mathbb{R}_{+}^{N}\}$ be an $\mathbb{R}_{+}^{N}-$parameter L\'{e}vy
process on $\mathbb{R}^{d}$. Then, according to O. E. Barndorff-Nielsen \textit{et al} \cite{Barndorff} Theorem 4.5, there exist N independent
L\'{e}vy processes $\{Z_{j}(t); t\geq0\}$ ($j=1,\cdots,N$) on $\mathbb{R}^{d}$ such that the $\mathbb{R}_{+}^{N}-$parameter L\'{e}vy process 
$\{V(\textbf{s}); \textbf{s}\in\mathbb{R}_{+}^{N}\}$ defined for $\textbf{s}=(s_{1},\cdots,s_{N})\in\mathbb{R}_{+}^{N}$ by 
$V(\textbf{s})=Z_{1}(s_{1})+ \cdots +Z_{N}(s_{N})$, has the property that, for any choice of $m\geq1$ and $\textbf{s}^{1}\preceq\cdots\preceq\textbf{s}^{m}$,
$$(X(\textbf{s}^{1}),\cdots,X(\textbf{s}^{m}))\stackrel{(law)}{=}(V(\textbf{s}^{1}),\cdots,V(\textbf{s}^{m})).$$
Now, for all $m,n\in\mathbb{N}^{\ast}$ and $\theta =(\theta_{1},\cdots,\theta_{m})\in\mathbb{R}^{m}$, we have:
$$\mathbb{E}exp\{i\sum_{k=1}^{m}<\theta_{k}, X(n.\textbf{s}^{k})>\}=\mathbb{E}exp\{i\sum_{k=1}^{m}<\theta_{k}, V(n.\textbf{s}^{k})>\}.$$
This implies that:
$$\mathbb{E}exp\{i\sum_{k=1}^{m}<\theta_{k}, X(n.\textbf{s}^{k})>\}=\mathbb{E}exp\{i\sum_{k=1}^{m}<\theta_{k}, \sum_{j=1}^{N} Z_{j}(n.\textbf{s}^{k}_{j})>\}.$$
According to the independence of the processes $\{Z_{j}(t); t\geq0\}$, $j=1,\cdots,N$, we have:
$$\mathbb{E}exp\{i\sum_{k=1}^{m}<\theta_{k}, X(n.\textbf{s}^{k})>\}=\prod_{j=1}^{N}\mathbb{E}exp\{i\sum_{k=1}^{m}<\theta_{k},Z_{j}(n.\textbf{s}^{k}_{j})>\}.$$
Thanks to the IDT property of the L\'{e}vy processes $\{Z_{j}(t); t\geq0\}$, $j=1,\cdots,N$, we get
$$\mathbb{E}exp\{i\sum_{k=1}^{m}<\theta_{k}, X(n.\textbf{s}^{k})>\}=\prod_{j=1}^{N}[\mathbb{E}exp\{i\sum_{k=1}^{m}<\theta_{k},Z_{j}(\textbf{s}^{k}_{j})>\}]^{n}.$$
That is,
$$\mathbb{E}exp\{i\sum_{k=1}^{m}<\theta_{k}, X(n.\textbf{s}^{k})>\}=[\prod_{j=1}^{N}\mathbb{E}exp\{i\sum_{k=1}^{m}<\theta_{k},Z_{j}(\textbf{s}^{k}_{j})>\}]^{n}.$$
Hence, by the independence of the L\'{e}vy processes $\{Z_{j}(t); t\geq0\}$, $j=1,\cdots,N$, we get:
$$\mathbb{E}exp\{i\sum_{k=1}^{m}<\theta_{k}, X(n.\textbf{s}^{k})>\}=[\mathbb{E}exp\{i\sum_{k=1}^{m}<\theta_{k}, \sum_{j=1}^{N} Z_{j}(\textbf{s}^{k}_{j})>\}]^{n}.$$
And then,
$$\mathbb{E}exp\{i\sum_{k=1}^{m}<\theta_{k}, X(n.\textbf{s}^{k})>\}=[\mathbb{E}exp\{i\sum_{k=1}^{m}<\theta_{k}, X(\textbf{s}^{k})>\}]^{n}.$$
This prove that X is a multiparameter IDT process of type 2. \\
$(\Longleftarrow)$
Let $\{X(\textbf{s});\textbf{s}=(s_{1},\cdots,s_{N})\in\mathbb{R}_{+}^{N}\}$ be a multiparameter IDT process of type 2 and assumed that it is 
stochastically continuous and provided the independence of increments. Then, it is enough to prove that for 
any $\textbf{s}^{1}\preceq\textbf{s}^{2}$ and $\textbf{s}^{3}\preceq\textbf{s}^{4}$ satisfying $\textbf{s}^{2}-\textbf{s}^{1}=\textbf{s}^{4}-\textbf{s}^{3}$, $$X(\textbf{s}^{2})-X(\textbf{s}^{1})\stackrel{(law)}{=}X(\textbf{s}^{4})-X(\textbf{s}^{3}).$$
Let $\textbf{s}\preceq\textbf{t}$, then $\textbf{t}-\textbf{s}=(t_{1}-s_{1},\cdots,t_{N}-s_{N})\in\mathbb{R}_{+}^{N}$ 
and $\mathcal{L}(X(\textbf{t}-\textbf{s}))=\mu_{\textbf{t}-\textbf{s}}$ is infinitely divisible and its triplet is such that:
$$A_{\textbf{t-s}}=(t_{1}-s_{1})A_{e^{1}}+ \cdots + (t_{N}-s_{N})A_{e^{N}}=A_{\textbf{t}}-A_{\textbf{s}},$$
$$\nu_{\textbf{t-s}}=(t_{1}-s_{1})\nu_{e^{1}}+ \cdots + (t_{N}-s_{N})\nu_{e^{N}}=\nu_{\textbf{t}}-\nu_{\textbf{s}},$$
$$\gamma_{\textbf{t-s}}=(t_{1}-s_{1})\gamma_{e^{1}}+ \cdots + (t_{N}-s_{N})\gamma_{e^{N}}=\gamma_{\textbf{t}}-\gamma_{\textbf{s}}.$$
And then, for all $z\in\mathbb{R}^{d}$, we have
$$\hat{\mu}_{\textbf{t-s}}(z)=exp\{-\frac{1}{2}<z,A_{\textbf{t-s}} z> + \int_{\mathbb{R}^{d}} g(z,x)\nu_{\textbf{t-s}}(dx)+i<z,\gamma_{\textbf{t-s}}> \},$$
with $g(z,x)=e^{i<z,x>}-1-i<z,x>1_{\{\mid x \mid \leq 1\}}(x)$. \\
This implies that,
$$\hat{\mu}_{\textbf{t-s}}(z)=[\hat{\mu}_{\textbf{t}}(z)][\hat{\mu}_{\textbf{s}}(z)]^{-1}.$$
That is,
$$\mathbb{E}exp\{i<z,X(\textbf{t}-\textbf{s})>\}=
[\mathbb{E}exp\{i<z,X(\textbf{t})-X(\textbf{s})+X(\textbf{s})>\}][\mathbb{E}exp\{i<z,X(\textbf{s})>\}]^{-1}.$$
Now, according to the assumption of independence of increments, we get:
$$\mathbb{E}exp\{i<z,X(\textbf{t}-\textbf{s})>\}=\mathbb{E}exp\{i<z,X(\textbf{t})-X(\textbf{s})>\}.$$
Hence, for any $\textbf{s}^{1}\preceq\textbf{s}^{2}$ and $\textbf{s}^{3}\preceq\textbf{s}^{4}$ satisfying $\textbf{s}^{2}-\textbf{s}^{1}=\textbf{s}^{4}-\textbf{s}^{3}$, we have $$X(\textbf{s}^{2})-X(\textbf{s}^{1})\stackrel{(law)}{=}X(\textbf{s}^{4})-X(\textbf{s}^{3}).$$
In other hand, one may easily point that $X(\textbf{0})=0$ almost surely. \\
In fact, for $n=2$, we get:
$$\hat{\mu}_{2.\textbf{0}}(z)=(\hat{\mu}_{\textbf{0}}(z))^{2} \Rightarrow \hat{\mu}_{\textbf{0}}(z)=1=\hat{\delta}_{0} \mbox{ i.e. } \mathcal{L}(X(\textbf{0}))=\delta_{0}.$$
The proof is now completed.

\end{proof}

\begin{rem}

The L\'{e}vy's $\mathbb{R}_{+}^{M}$-parameter Brownian motion emphasizes more again on Theorem 8.1. In fact, J. Pedersen and K. Sato showed in
$\cite{PedersenSato}$ Example 2.16, that the L\'{e}vy's $\mathbb{R}_{+}^{M}$-parameter Brownian motion is not an $\mathbb{R}_{+}^{M}$-parameter
L\'{e}vy process for default of the independence of increments.

\end{rem}

\begin{thm}

Let $\{X(\textbf{s});\textbf{s}=(s_{1},\cdots,s_{N})\in\mathbb{R}_{+}^{N}\}$ be a stochastically continuous multiparameter IDT process of type 2.
Then, there exist an associated $\mathbb{R}_{+}^{N}-$parameter L\'{e}vy process $\{L(\textbf{s});\textbf{s}\in\mathbb{R}_{+}^{N}\}$ i.e. for any fixed $\textbf{s}\in\mathbb{R}_{+}^{N}$,
$$\mathcal{L}(X(\textbf{s}))=\mathcal{L}(X(s_{1},\cdots,s_{N}))=\mathcal{L}(L(s_{1},\cdots,s_{N}))
=\mathcal{L}(L(\textbf{s})).$$

\end{thm}

\begin{proof}

By Proposition \ref{prop8.1}, we have that $\{\mu(\textbf{s}); \textbf{s}\in\mathbb{R}_{+}^{N}\}$ is an $\mathbb{R}_{+}^{N}-$parameter convolution semigroup
on $\mathbb{R}^{d}$. For $j=1,\cdots,N$, $e^{j}=(\delta_{jk})_{1\leq k\leq N}$, where $\delta_{jk}=0$ or 1 according as $k \neq j$ or $k=j$,
form a strong basis of $\mathbb{R}_{+}^{N}$. Then, it follows by Theorem 3.2 in J. Pedersen and K. Sato $\cite{PedersenSato}$,
that $\{\mu(\textbf{s}); \textbf{s}\in\mathbb{R}_{+}^{N}\}$ is generative i.e. there exist an $\mathbb{R}_{+}^{N}-$parameter L\'{e}vy process
(in law) $L=\{L(\textbf{s}); \textbf{s}\in\mathbb{R}_{+}^{N}\}$ such that for any $\textbf{s}\in\mathbb{R}_{+}^{N}$, we have
$\mathcal{L}(L(\textbf{s}))=\mu (\textbf{s})=\mathcal{L}(X(\textbf{s}))$. The proof is complete.

\end{proof}

\begin{rem}

We prove that Brownian sheet studied by S. Orey and W. E. Pruitt \cite{Orey}, M. Talagrand \cite{Talagrand}, and many others, and the so-called 
multiparameter L\'{e}vy process studied by W. Ehm \cite{Ehm} (in the strictly stable case), M. E. Vares \cite{Vares}, and S. Lagaise \cite{Lagaise} (both in the case N=2), are not multiparameter IDT of type 2. \\
In fact, if $\{X(\textbf{s}); \textbf{s}\in\mathbb{R}_{+}^{2}\}$ is a L\'{e}vy sheet, $\mathcal{L}(X(\textbf{s}))$ is
infinitely divisible, and denoting $\mu$ the law of $X(1,1)$, we have
$$\mathbb{E}exp\{i<z,X(\textbf{s})>\}=\hat{\mu}(z)^{s_{1}s_{2}} \mbox{ for } \textbf{s}=(s_{1},s_{2})\in\mathbb{R}_{+}^{2}.$$
If $\{X(\textbf{s}); \textbf{s}\in\mathbb{R}_{+}^{2}\}$ is a bi-parameter IDT process of type 2, then $\{\mathcal{L}(X(\textbf{s})); \textbf{s}\in\mathbb{R}_{+}^{2}\}$ is an $\mathbb{R}_{+}^{2}$-parameter convolution semigroup.
Following J. Pedersen and K. Sato \cite{PedersenSato} Theorem 1.2,
$$\mathbb{E}exp\{i<z,X(\textbf{s})>\}=[\mathbb{E}exp\{i<z,X(1,0)>\}]^{s_{1}}[\mathbb{E}exp\{i<z,X(0,1)>\}]^{s_{2}}$$
and this is different from $\hat{\mu}(z)^{s_{1}s_{2}}$ for $s_{1}\neq s_{2}$.

\end{rem}

\subsection{Links with operator stability and operator selfsimilarity}

In order to connect multiparameter IDT processes of type 2 with operator stability and operator selfsimilarity, we recall definitions of these concepts (in the sense of M. Maejima \cite{Maejima1995}, and K. Sato \cite{Sato1991} and \cite{Sato2004}).

\begin{defn}

An $\mathbb{R}^{d}$-valued multiparameter process $\{X(\textbf{s}); \textbf{s}\in\mathbb{R}_{+}^{N}\}$ is operator selfsimilar with exponent Q or 
Q-selfsimilar if and only if for every $\alpha>0$,
$$\{X(\alpha\textbf{s}); \textbf{s}\in\mathbb{R}_{+}^{N}\}\stackrel{(law)}{=}
\{{\alpha}^{Q}X(\textbf{s}); \textbf{s}\in\mathbb{R}_{+}^{N}\}.$$

\end{defn}

\begin{defn}

An $\mathbb{R}^{d}$-valued multiparameter process $\{X(\textbf{s}); \textbf{s}\in\mathbb{R}_{+}^{N}\}$ is strictly operator stable with exponent Q or strictly Q-stable if and only if for every positive integer n,
$$\{\sum_{i=1}^{n} X^{(i)}(\textbf{s}); \textbf{s}\in\mathbb{R}_{+}^{N}\}\stackrel{(law)}{=}\{n^{Q}X(\textbf{s}) ;\textbf{s}\in\mathbb{R}_{+}^{N} \},$$
where $\{X^{(i)}(\textbf{s}); \textbf{s}\in\mathbb{R}_{+}^{N}\}$, $i=1,\cdots,n$, are independent copies of $\{X(\textbf{s}); \textbf{s}\in\mathbb{R}_{+}^{N}\}$.

\end{defn}

Now, we have the following extensions of Proposition 3.1 and Theorem 3.2 in K. Es-Sebaiy and Y. Ouknine $\cite{Ouknine}$.

\begin{prop}

A nontrivial, strictly Q-stable, Q-selfsimilar multiparameter process $X=\{X(\textbf{s});\textbf{s}=(s_{1},\cdots,s_{N})\in\mathbb{R}_{+}^{N}\}$ is a
multiparameter IDT of type 2.

\end{prop}

\begin{proof}

Since X is strictly Q-stable, we have for all $n\in\mathbb{N}^{\ast}$,
$$\{n^{Q} X(\textbf{s}); \textbf{s}\in\mathbb{R}_{+}^{N}\}\stackrel{(law)}{=}\{\sum_{i=1}^{n} X^{(i)}(\textbf{s}); \textbf{s}\in\mathbb{R}_{+}^{N}\},$$
where $X^{(i)}$ for $i=1,\cdots,n$, are independent copies of X. \\
Now, according to the Q-selfsimilarity of X, we get:
$$\{n^{Q} X(\textbf{s}); \textbf{s}\in\mathbb{R}_{+}^{N}\}\stackrel{(law)}{=}\{X(n \textbf{s}); \textbf{s}\in\mathbb{R}_{+}^{N}\}.$$
This implies that X is a multiparameter IDT process of type 2, i.e. for all $n\in\mathbb{N}^{\ast}$
$$\{X(n \textbf{s}); \textbf{s}\in\mathbb{R}_{+}^{N}\}\stackrel{(law)}{=}\{\sum_{i=1}^{n} X^{(i)}(\textbf{s}); \textbf{s}\in\mathbb{R}_{+}^{N}\},$$
where $X^{(i)}$ for $i=1,\cdots,n$, are independent copies of X.

\end{proof}

\begin{thm}

Let $\{X(\textbf{s}); \textbf{s}=(s_{1},\cdots,s_{N})\in\mathbb{R}_{+}^{N}\}$ be a non trivial, stochastically continuous multiparameter IDT process of type 2. Then, $\{X(\textbf{s}); \textbf{s}\in\mathbb{R}_{+}^{N}\}$ is strictly Q-stable if and only if it is Q-selfsimilar.

\end{thm}

\begin{proof}

Let $X=\{X(\textbf{s}); \textbf{s}\in\mathbb{R}_{+}^{N}\}$ be a stochastically continuous multiparameter IDT process of type 2. \\
$(\Longleftarrow)$ Assume that X is Q-selfsimilar. Then for all $n\in\mathbb{N}^{\ast}$,
$$\{X(n \textbf{s}); \textbf{s}\in\mathbb{R}_{+}^{N}\}\stackrel{(law)}{=}\{n^{Q} X(\textbf{s}); \textbf{s}\in\mathbb{R}_{+}^{N}\}.$$
Now, since X is a multiparameter IDT process of type 2, we get:
$$\{n^{Q} X(\textbf{s}); \textbf{s}\in\mathbb{R}_{+}^{N}\}\stackrel{(law)}{=}\{X(n \textbf{s}); \textbf{s}\in\mathbb{R}_{+}^{N}\}
\stackrel{(law)}{=}\{\sum_{i=1}^{n} X^{(i)}(\textbf{s}); \textbf{s}\in\mathbb{R}_{+}^{N}\},$$
where $X^{(i)}$ for $i=1,\cdots,n$, are independent copies of X.\\
Then X is strictly Q-stable. \\
$(\Longrightarrow)$ Assume that X is strictly Q-stable. \\
Then, for all $n\in\mathbb{N}^{\ast}$, we have:
$$\{n^{Q} X(\textbf{s}); \textbf{s}\in\mathbb{R}_{+}^{N}\}\stackrel{(law)}{=}\{\sum_{i=1}^{n} X^{(i)}(\textbf{s}); \textbf{s}\in\mathbb{R}_{+}^{N}\},$$
where $X^{(i)}$ for $i=1,\cdots,n$, are independent copies of X. \\
Since X is a multiparameter IDT process of type 2, it follows that for all $n\in\mathbb{N}^{\ast}$,
$$\{n^{Q} X(\textbf{s}); \textbf{s}\in\mathbb{R}_{+}^{N}\}\stackrel{(law)}{=}\{X(n \textbf{s}); \textbf{s}\in\mathbb{R}_{+}^{N}\}.$$
An then, we get:
$$\{(\frac{1}{n})^{Q} X(\textbf{s}); \textbf{s}\in\mathbb{R}_{+}^{N}\}\stackrel{(law)}{=}\{X(\frac{1}{n}\textbf{s}); \textbf{s}\in\mathbb{R}_{+}^{N}\}.$$
Hence, for any $m,n\in\mathbb{N}^{\ast}$, we have
$$\{(\frac{m}{n})^{Q} X(\textbf{s}); \textbf{s}\in\mathbb{R}_{+}^{N}\}\stackrel{(law)}{=}\{X(\frac{m}{n}\textbf{s}); \textbf{s}\in\mathbb{R}_{+}^{N}\}.$$
Now, thanks to the stochastically continuity of X and the density of $\mathbb{Q}_{+}$ in $\mathbb{R}_{+}$, we obtain for all $\alpha >0$,
$$\{{\alpha}^{Q} X(\textbf{s}); \textbf{s}\in\mathbb{R}_{+}^{N}\}\stackrel{(law)}{=}\{X(\alpha \textbf{s}); \textbf{s}\in\mathbb{R}_{+}^{N}\}.$$
Then X is Q-selfsimilar and this complete the proof.

\end{proof}

\subsection{Links with temporal selfdecomposability and operator selfdecomposability}

In this section, we are interested to link multiparameter IDT processes of type 2 with temporal selfdecomposability and operator selfdecomposability. First, we recall the definitions of these concepts.

\begin{defn}

An $\mathbb{R}^{d}$-valued multiparameter process $X=\{X(\textbf{s}); \textbf{s}\in\mathbb{R}_{+}^{N}\}$ is temporally selfdecomposable of order 1, if and only if, for any $c\in (0,1)$,
there exist a process $U^{(c)}=\{U^{(c)}(\textbf{s}); \textbf{s}\in\mathbb{R}_{+}^{N}\}$ (called the c-residual of X) independent of $X$ and such that
\begin{equation}\label{17}
\{X(\textbf{s}); \textbf{s}\in\mathbb{R}_{+}^{N}\}\stackrel{(law)}{=}\{X(c \textbf{s})+U^{(c)}(\textbf{s}); \textbf{s}\in\mathbb{R}_{+}^{N}\}.
\end{equation}

For every $m>1$, X is temporally selfdecomposable of order m, if for any $c\in (0,1)$, the c-residual $U^{(c)}$ is temporally selfdecomposable of order $m-1$. When X is temporally selfdecomposable of order m for all $m\in\mathbb{N}^{\ast}$, X is said temporally selfdecomposable of infinite order.

\end{defn}

\begin{defn}\label{def8.5}

An $\mathbb{R}^{d}$-valued multiparameter process $X=\{X(\textbf{s}); \textbf{s}\in\mathbb{R}_{+}^{N}\}$ is Q-selfdecomposable or Q operator selfdecomposable, if and only if for every $c\in (0,1)$,
\begin{equation}\label{18}
X \stackrel{(law)}{=} c^{Q} X'+U^{(c)},
\end{equation}
where $X'=\{X'(\textbf{s}); \textbf{s}\in\mathbb{R}_{+}^{N}\}$ is a copie of X, $U^{(c)}=\{U^{(c)}(\textbf{s}); \textbf{s}\in\mathbb{R}_{+}^{N}\}$ a multiparameter stochastic process on $\mathbb{R}^{d}$, and X' and $U^{(c)}$ are independent.

\end{defn}

\begin{prop}\label{prop8.4}

An $\mathbb{R}^{d}$-valued stochastically continuous multiparameter IDT process of type 2, is temporally selfdecomposable of infinite order.

\end{prop}

\begin{proof}

Let $m\in\mathbb{N}^{\ast}$ and $X=\{X(\textbf{s}); \textbf{s}\in\mathbb{R}_{+}^{N}\}$ a multiparameter IDT process of type 2.
Assume that X is stochastically continuous, then for any $\textbf{s}^{1},\cdots,\textbf{s}^{m}\in\mathbb{R}_{+}^{N}$, and for all
$\theta_{1},\cdots,\theta_{m}\in\mathbb{R}$ and $c\in (0,1)$, we get
$$\mathbb{E} exp\{i\sum_{k=1}^{m} <\theta_{k}, X(\textbf{s}^{k})>\}=(\mathbb{E} exp\{i\sum_{k=1}^{m} <\theta_{k}, X( c \textbf{s}^{k})>\})^{1/c}.$$
This implies that:
$$\mathbb{E} exp\{i\sum_{k=1}^{m} <\theta_{k}, X(\textbf{s}^{k})>\}=(\mathbb{E} exp\{i\sum_{k=1}^{m} <\theta_{k}, X( c \textbf{s}^{k})>\})^{1/c'}
(\mathbb{E} exp\{i\sum_{k=1}^{m} <\theta_{k}, X(c \textbf{s}^{k})>\}),$$
where $\frac{1}{c'}=\frac{1}{c}-1$. \\
Then, it follows that:
$$\mathbb{E} exp\{i\sum_{k=1}^{m} <\theta_{k}, X(\textbf{s}^{k})>\}=(\mathbb{E} exp\{i\sum_{k=1}^{m} <\theta_{k}, X( c \textbf{s}^{k})>\})
(\mathbb{E} exp\{i\sum_{k=1}^{m} <\theta_{k}, X(\frac{c}{c'}\textbf{s}^{k})>\}).$$
Therefore X is temporally selfdecomposable and
$$\{X(\textbf{s}); \textbf{s}\in\mathbb{R}_{+}^{N}\}\stackrel{(law)}{=}\{X(c \textbf{s}) + U^{(c)}(\textbf{s}); \textbf{s}\in\mathbb{R}_{+}^{N}\},$$
where $$\{U^{(c)}(\textbf{s}); \textbf{s}\in\mathbb{R}_{+}^{N}\}\stackrel{(law)}{=}\{ X(\frac{c}{c'} \textbf{s}); \textbf{s}\in\mathbb{R}_{+}^{N}\}.$$
It is clear that $U^{(c)}$ is a stochastically continuous multiparameter IDT process of type 2. The same steps as above applied to $U^{(c)}$, proves that
it is temporally selfdecomposable and so on, the desired result follows.

\end{proof}

\begin{prop}\label{prop8.5}

Let $X=\{X(\textbf{s}); \textbf{s}\in\mathbb{R}_{+}^{N}\}$ an $\mathbb{R}^{d}$-valued multiparameter IDT process of type 2.
Then, if X is Q-selfdecomposable, the process $U^{(c)}$ defined in (\ref{18}) is a multiparameter IDT process of type 2.

\end{prop}

\begin{proof}

Let $m,n\in\mathbb{N}^{\ast}$ and $\theta_{k}\in\mathbb{R}$ for $k=1,\cdots,m$. \\
According to Definition \ref{def8.5} above, we get:
$$\mathbb{E} exp\{i\sum_{k=1}^{m} <\theta_{k}, X(\textbf{s}^{k})>\}=\mathbb{E} exp\{i\sum_{k=1}^{m} <\theta_{k}, c^{Q} X'(\textbf{s}^{k})+U^{(c)}(\textbf{s}^{k})>\}.$$
The independence of $U^{(c)}$ and $X'$ implies that:
$$\mathbb{E} exp\{i\sum_{k=1}^{m} <\theta_{k}, X(\textbf{s}^{k})>\}=\mathbb{E} exp\{i\sum_{k=1}^{m} <\theta_{k}, c^{Q} X'(\textbf{s}^{k})>\}
\mathbb{E} exp\{i\sum_{k=1}^{m} <\theta_{k}, U^{(c)}(\textbf{s}^{k})>\}.$$
And then we get,
$$\mathbb{E} exp\{i\sum_{k=1}^{m} <\theta_{k}, U^{(c)}(\textbf{s}^{k})>\}=\frac{\mathbb{E} exp\{i\sum_{k=1}^{m} <\theta_{k}, X(\textbf{s}^{k})>\}}
{\mathbb{E} exp\{i\sum_{k=1}^{m} <\theta_{k}, c^{Q} X'(\textbf{s}^{k})>\}}.$$
Since X and X' are multiparameter IDT processes of type 2, it follows that
$$\mathbb{E} exp\{i\sum_{k=1}^{m} <\theta_{k}, U^{(c)}(n \textbf{s}^{k})>\}=\frac{(\mathbb{E} exp\{i\sum_{k=1}^{m}
<\theta_{k}, X(\textbf{s}^{k})>\})^{n}}{(\mathbb{E} exp\{i\sum_{k=1}^{m} <\theta_{k}, c^{Q} X'(\textbf{s}^{k})>\})^{n}}.$$
Then, $U^{(c)}$ is a multiparameter IDT process of type 2, that is:
$$\mathbb{E} exp\{i\sum_{k=1}^{m} <\theta_{k}, U^{(c)}(n \textbf{s}^{k})>\}=
(\mathbb{E} exp\{i\sum_{k=1}^{m} <\theta_{k}, U^{(c)}(\textbf{s}^{k})>\})^{n}.$$

\end{proof}

\begin{coro}

An $\mathbb{R}^{d}$-valued multiparameter IDT process of type 2 which is continuous in probability and Q-selfsimilar, is Q-selfdecomposable.

\end{coro}

\begin{proof}

The corollary is a consequence of Propositions \ref{prop8.4} and \ref{prop8.5}.

\end{proof}

\subsection{Subordination by a multiparameter IDT process of type 2}

The aim of this section, is to investigate the subordination of an $\mathbb{R}_{+}^{N}-$parameter L\'{e}vy process by a chronometer
(i.e. an increasing and stochastically continuous process starting at 0), which is a multiparameter IDT process of type 2.

\begin{prop}

Let $\{X(\textbf{s}); \textbf{s}=(s_{1},\cdots,s_{N})\in\mathbb{R}_{+}^{N}\}$ be an $\mathbb{R}_{+}^{N}-$parameter L\'{e}vy process
on $\mathbb{R}^{d}$ and consider 
$\{\xi(\textbf{s});\textbf{s}=(s_{1},\cdots,s_{N})\}=
\{(\xi^{1}_{s_{1}},\cdots,\xi^{N}_{s_{N}})\}$ a multiparameter IDT process of type 2
with $\{\xi^{j}_{t}; t\in\mathbb{R}_{+}\}$ are independent IDT chronometers. Assume that $\{\xi(\textbf{s})\}$ is independent of $\{X(\textbf{s})\}$
and define the subordinated process by composition as follows:
$$Y(\textbf{s})=X(\xi_{\textbf{s}}),\,\,\textbf{s}\in\mathbb{R}_{+}^{N}.$$
Then, $\{Y(\textbf{s}); \textbf{s}\in\mathbb{R}_{+}^{N}\}$ is a multiparameter IDT process of type 2.

\end{prop}

\begin{proof}

Let $\xi^{(l)}$, $l=1,\cdots,n$, be independent copies of $\xi$. Since X is independent of $\xi$, then for every $m\geq1$ and
$\theta=(\theta_{1},\cdots,\theta_{m})\in\mathbb{R}^{m},$ we have
$$J(n,\theta):=\mathbb{E}exp\{i\sum_{k=1}^{m} <\theta_{k}, X(\xi_{n\textbf{t}^{k}})>\}
=\mathbb{E}[(\mathbb{E}exp \{i\sum_{k=1}^{m} <\theta_{k}, X(\textbf{s}^{k})>\})_{\textbf{s}^{k}=\xi_{n\textbf{t}^{k}},\, k=1,\cdots,m}].$$
Using the IDT property, we obtain
$$J(n,\theta)=\mathbb{E}\big[(\mathbb{E}exp \{i\sum_{k=1}^{m}
<\theta_{k}, X(\textbf{s}^{k})>\})_{\textbf{s}^{k}=\sum_{l=1}^{n}\xi^{(l)}_{\textbf{t}^{k}},\, k=1,\cdots,m}\big].$$
According to the change of variables $\lambda_{k}=\theta_{k} + \cdots + \theta_{m}$ and $\textbf{t}^{0}=\textbf{0}$,  we have
$$J(n,\theta)=\mathbb{E}[(\mathbb{E}exp \{i\sum_{k=1}^{m} <\lambda_{k}, X(\textbf{s}^{k})-X(\textbf{s}^{k-1})>\})_{\textbf{s}^{k}
=\sum_{l=1}^{n}\xi^{(l)}_{\textbf{t}^{k}},\, k=1,\cdots,m}].$$
By the independence of increments of X, we get
$$J(n,\theta)=\mathbb{E}[(\prod_{k=1}^{m} \mathbb{E}exp
\{i <\lambda_{k}, X(\textbf{s}^{k})-X(\textbf{s}^{k-1})>\})_{\textbf{s}^{k}=\sum_{l=1}^{n}\xi^{(l)}_{\textbf{t}^{k}}, \, k=1,\cdots,m}].$$
Now, it follows from the stationary of the increments of X and the independence of the $\xi^{(l)}$, $l=1,\cdots,n$, that
$$J(n,\theta)=\mathbb{E}[(\prod_{k=1}^{m}\prod_{l=1}^{n} \mathbb{E}exp
\{i <\lambda_{k}, X(\textbf{r}^{l,k})>\})_{\textbf{r}^{l,k}=\xi^{(l)}_{\textbf{t}^{k}}-\xi^{(l)}_{\textbf{t}^{k-1}},\,
k=1,\cdots,m;\, l=1,\cdots,n}].$$
And that is
$$J(n,\theta)=\mathbb{E}[(\prod_{l=1}^{n} \mathbb{E}exp \{i \sum_{k=1}^{m}
<\theta_{k}, X(\textbf{r}^{l,k})>\})_{\textbf{r}^{l,k}=\xi^{(l)}_{\textbf{t}^{k}},\, k=1,\cdots,m;\, l=1,\cdots,n}].$$
Now, this implies
$$J(n,\theta)=(\mathbb{E}[(\mathbb{E}exp \{i \sum_{k=1}^{m} <\theta_{k}, X(\xi_{\textbf{r}^{k}})>\})_{\textbf{r}^{k}
=\xi_{\textbf{t}^{k}},\, k=1,\cdots,m}])^{n}.$$
Hence
$$J(n,\theta)=(\mathbb{E}exp \{i \sum_{k=1}^{m} <\theta_{k}, X(\xi_{\textbf{t}^{k}})>\})^{n}.$$
The proof is achieved.

\end{proof}

\textit{\textbf{Acknowledgements.}}
The authors are grateful to professor B. Roynette for his various and pertinent suggestions.
We also thank him for the graceful monograph \textit{Peacocks and associated martingales with explicit construction},
which already offers a fertile field for a fruitful development of IDT processes. We also thanks the anonymous referee for
careful reading and significative remarks which allowed us to improve the paper.


\begin{thebibliography}{1}

\bibitem{Adler}
R. J. Adler, D. Monrad, R. H. Scissors and R. Wilson,
\emph{Representations, decompositions and sample function continuity of random fields with independent increments, Stochastic Processes and their Applications, 1983, 3-30}

\bibitem{Applebaum}
D. Applebaum,
\emph{Lévy processes and stochastic calculus, Cambridge studies in advanced mathematics, 2004, Cambridge University Press}

\bibitem{K.Sato}
O. E. Barndorff-Nielsen, M. Maejima and K. Sato,
\emph{Infinite Divisibility for Stochastic Processes and Time Change, Journal of Theorical Probability, 2006, 19, 411-446}

\bibitem{Barndorff}
O. E. Barndorff-Nielsen, J. Pedersen and K. Sato,
\emph{Multivariate subordination, selfdecomposability and stability, Adv. Appl. Probab. 2001, 160-187}

\bibitem{Bochner}
S. Bochner,
\emph{Diffusion equation and stochastic processes, P. Nat. Acad. Sci., 1949, 85, 369-370}

\bibitem{Chentsov}
N. N. Chentsov,
\emph{Lévy's Brownian motion of several parameters and generalized white noise, Theory Probab. Appl. 1957, 265-266}

\bibitem{Dalang}
R. C. Dalang and J. B. Walsh,
\emph{The sharp Markov property of Lévy sheets, Annals of Probability, 1992, 20, 2, 591-626}

\bibitem{Ehm}
W. Ehm, \emph{Sample function properties of multi-parameter stable processes, Z. Wahrsch. Verw. Gebiete, 1981, 499-530}

\bibitem{embrechts2002selfsimilar}
P. Embrechts and M. Maejima,
\emph{Selfsimilar processes, Princeton series in applied mathematics, 2002, Princeton University Press}

\bibitem{Ouknine}
K. Es-Sebaiy and Y. Ouknine,
\emph{How rich is the class of processes which are infinitely divisible with respect to time ? 
Statistics and Probability Letters, 2008, 78, 537-547}

\bibitem{Follmer}
H. Föllmer, C.T. Wu and M. Yor,
\emph{On weak Brownian motions of arbitrary order, Annales de l'Institut Henri Poincaré (B) Probability and Statistics, 2000, 36, 447-487}

\bibitem{Hirsch}
F. Hirsch, B. Roynette and M. Yor,
\emph{Unifying constructions of martingales associated with processes increasing in the convex order, via Lévy and Sato sheets. 
Expositiones Mathematicae, 2010, 28, 4, 299-324}

\bibitem{Roynette}
F. Hirsch, B. Roynette,and M. Yor,
\emph{From an Itô type calculus for Gaussian processes to integrals of log-normal processes increasing in the convex order, J. Mat. Soc. Japan, 2010}

\bibitem{M.Yor}
F. Hirsch, B. Roynette and M. Yor,
\emph{Applying Itô's motto: "look at the infinite dimensional picture" by constructing sheets to obtain processes increasing in the convex order,
Periodica Mathematica Hungarica, 2010, 61, 195-211}

\bibitem{Yor}
F. Hirsch and M. Yor,
\emph{A construction of processes with one dimensional martingale marginals, based upon path-space Ornstein-Uhlenbeck processes and the Brownian sheet,
J. Math. Kyoto Univ., 2009, 49, 2, 389-417}

\bibitem{HirschYor}
F. Hirsch and M. Yor,
\emph{A construction of processes with one dimensional martingale marginals, associated with a Lévy process, via it's Lévy sheet, J. Math. Kyoto Univ.,
2009, 49, 4, 785-815}

\bibitem{Hirschbook}
F. Hirsch , C. Profeta , B. Roynette  and M. Yor,
\emph{Peacocks and Associated Martingales, with Explicit Constructions, Bocconi and Springer Series, 2011, Springer}

\bibitem{Lagaise}
S. Lagaize, \emph{Hölder exponent for a two-parameter Lévy process, J. Multivariate Anal., 2001, 270-285}

\bibitem{Levy}
P. Lévy, \emph{Processus Stochastiques et Mouvement Brownien, Gauthier-Villars, Paris, 1948 (2e éd. 1965)}

\bibitem{Maejima}
M. Maejima and K. Sato,  \emph{Semi-selfsimilar processes, J. Th. Prob. 1999, 347-373}

\bibitem{Maejima1995}
M. Maejima,  \emph{Operator stable processes and operator fractional stable processes, Probability and Mathematical Statistics, 1995, 449-460}

\bibitem{Mansuy}
R. Mansuy,  \emph{On processes which are infinitely divisible with respect to time, arXiv:math/0504408, 2005}

\bibitem{Mckean}
H. P. McKean Jr.,  \emph{Brownian motion with a several dimensional time, Theory Probab. Appl., 1963, 357-378}

\bibitem{Orey}
S. Orey and W. E. Pruitt,  \emph{Samples functions of the N-parameters Wiener process, Ann. Probab., 1973, 138-163}

\bibitem{SatoPedersen}
J. Pedersen and K. Sato, \emph{Cone-parameter convolution semigroups and their subordination, Tokyo J. Math., 2003, 503-525}

\bibitem{PedersenSato}
J. Pedersen and K. Sato, \emph{Semigroups and processes with parameter in a cone, Abstract and Applied Analysis (ed. N. M. Chuong et al., World Scientific), 2004, 499-513}

\bibitem{Pedersen}
J. Pedersen and K. Sato, \emph{Relation between cone-parameter Lévy processes and convolution semigroups, J. Math. Soc. Japan, 2004, 541-559}

\bibitem{Protter}
P. E. Protter, \emph{Stochastic Integration and Differential Equations, Applications of Mathematics, Springer, 2004}

\bibitem{Rohatgi}
V. Rohatgi, F. Steutel and G. Székely, \emph{Infinite divisibility of products and quotients of i.i.d. random variables, Math. Sci., 1990, 53-59}

\bibitem{Sato}
K. Sato, \emph{Lévy processes and infinitely divisible distributions, Cambridge studies in advanced mathematics, 1999, Cambridge University Press}

\bibitem{Sato1991}
K. Sato, \emph{Self-similar processes with independent increments, Prob. Th. Rel. Fields, 1991, 285-300}

\bibitem{Sato2004}
K. Sato, \emph{Stochastic integrals in additive processes and application to semi-Lévy processes, Osaka J. Math., 2004, 41, 211-236}

\bibitem{Shanbhag}
D. N. Shanbhag, D. Pestana and M. Sreehari,
\emph{Some further results in infinite divisibility, Math. Proc. Camb. Phil. Soc., 1977, 289-295}

\bibitem{Talagrand}
M. Talagrand, \emph{The small ball problem for the Brownian sheet, Ann. Prob., 1994, 1331-1354}

\bibitem{Vares}
M. E. Vares, \emph{Local times for two-parameter Lévy processes, Stochastic process. Appl., 1983, 59-82}


\end{thebibliography}
\end{document}